\documentclass{gtart}

\def\ifplaintex{\expandafter\ifx\csname documentclass\endcsname\relax}

\ifplaintex 
\else
\fi

\expandafter\ifx\csname beginpicture\endcsname\relax
\expandafter\ifx\csname documentclass\endcsname\relax
\input pictex \else
\input prepictex \input pictex \input postpictex \fi\fi

\def\gt{{\mathsurround=0pt\it $\cal G\mskip-2mu$eometry \&\ 
$\cal T\!\!$opology}}        

\def\gtp{{\mathsurround=0pt\it $\cal G\mskip-2mu$eometry \&\ 
$\cal T\!\!$opology $\cal P\!$ublications}}  

\def\lognumber#1{\def\thelognumber{#1}}
\def\volumenumber#1{\def\thevolumenumber{#1}}
\def\papernumber#1{\def\thepapernumber{#1}}
\def\volumeyear#1{\def\thevolumeyear{#1}}

\def\pagenumbers#1#2{\def\startpage{#1}\def\finishpage{#2}}
\def\published#1{\def\publishdate{#1}}
\def\proposed#1{\def\theproposer{#1}}
\def\seconded#1{\def\theseconders{#1}}
\def\received#1{\def\receiveddate{#1}}

\def\accepted#1{\def\accepteddate{#1}}

\def\coverauthors#1{\def\thecoverauthors{#1}}
\def\asciiauthors#1{\def\theasciiauthors{#1}}
\def\asciiaddress#1{\def\theasciiaddress{#1}}

\let\\\par\let\thelognumber\relax
\let\thevolumenumber\relax\let\thepapernumber\relax
\let\thevolumeyear\relax\let\thesamplenumber\relax\let\startpage\relax
\let\finishpage\relax\let\publishdate\relax\let\receiveddate\relax
\let\reviseddate\relax\let\accepteddate\relax\let\theasciititle\relax
\let\theasciiauthors\relax\let\theasciiaddress\relax
\let\theasciiabstract\relax
\let\theasciiemail\relax\let\theshortauthors\relax\let\theshorttitle\relax
\let\thecoverauthors\relax

\long\def\maketitlep{   

\count0=\startpage

\gt\hfill      
\beginpicture
\setcoordinatesystem units <0.33truein, 0.33truein> point at 2.2 0.9
\setplotsymbol ({$\cal G$})
\plotsymbolspacing=9truept
\circulararc 315 degrees from 0 1 center at 0 0
\setplotsymbol ({$\cal T$})
\circulararc 315 degrees from 1 -1 center at 1 0
\endpicture
\break
{\small\ifx\thesamplenumber\relax 
Volume \else Sample
\fi\thevolumenumber\ (\thevolumeyear)
\startpage--\finishpage\nl
Published: \publishdate}
\vglue 0.5truein plus 0.4fil minus 0.1truein

{\parskip=0pt\leftskip 0pt plus 1fil\def\\{\par\smallskip}{\ifplaintex\large
\else\Large\fi\bf\thetitle}\par\medskip}   

\vglue 0pt plus 0.1fil 

{\parskip=0pt\leftskip 0pt plus 1fil\def\\{\par}{\sc\theauthors}
\par\medskip}

\vglue 0pt plus 0.1fil 

{\small\parskip=0pt\let\newline\\
{\leftskip 0pt plus 1fil\def\\{\par}{\sl\theaddress}\par}
\expandafter\ifx\theemail\relax    
\relax\else\vglue 5pt plus 0.02fil minus 2pt\def\\{\stdspace{\rm 
and}\stdspace} 
\cl{Email:\stdspace\tt\theemail}\fi
\ifx\theurl\relax                  
\relax\else\vglue 5pt plus 0.02fil minus 2pt\def\\{\stdspace{\rm 
and}\stdspace}
\cl{URL:\stdspace\tt\theurl}\fi\par}

\vglue 7pt plus 0.3fil minus 3pt

{\bf Abstract}
\vglue 5pt plus 0.1fil minus 2pt

\theabstract

\vglue 7pt plus 0.3fil minus 3pt

{\bf AMS Classification numbers}\quad Primary:\quad \theprimaryclass

Secondary:\quad \thesecondaryclass

\vglue 5pt plus 0.3fil minus 2pt

{\bf Keywords}\quad \thekeywords

\vglue 10pt plus 0.5fil minus 5pt

{\small  Proposed: \theproposer\hfill Received: \receiveddate\nl
Seconded: \theseconders\hfill 
\ifx\reviseddate\relax                         
Accepted: \accepteddate                        
\else
Revised: \reviseddate                          
\fi}
\eject
}       

\let\maketitlepage\maketitlep
\let\maketitle\maketitlepage

\font\phead=cmsl9 scaled 950
\font=cmsl9 scaled 1050
\font\pnum=cmbx10 scaled 913
\font\lnum=cmbx10 
\font\pfoot=cmsl9 scaled 950
\font=cmsl9 scaled 1050
\ifplaintex
\headline{\vbox to 0pt{\vskip -4.5mm\line{\small\phead\ifnum
\count0=\startpage ISSN 1364-0380 (on line)
1465-3060 (printed) \hfill {\pnum\folio}\else\ifodd\count0\def\\{ }%
\ifx\theshorttitle\relax\thetitle\else\theshorttitle\fi\hfill{\pnum\folio}
\else\def\\{ and }{\pnum\folio}\hfill\ifx\theshortauthors\relax\theauthors
\else\theshortauthors\fi\fi\fi}\vss}}
\footline{\vbox to 0pt{\vglue 0mm\line{\small\pfoot\ifnum\count0=\startpage
\copyright\ \gtp\hfill\else
\gt, Volume \thevolumenumber\ (\thevolumeyear)\hfill\fi}\vss
}}
\else
\makeatletter
\def\@oddhead{{\small\lhead\ifnum\count0=\startpage ISSN 1364-0380 (on line)
1465-3060 (printed) \hfill {\lnum\number\count0}\else\ifodd\count0
\def\\{ }\ifx\theshorttitle\relax \thetitle \else\theshorttitle\fi\hfill
{\lnum\number\count0}\else\def\\{ and }{\lnum\number\count0}
\hfill\ifx\theshortauthors\relax 
\theauthors\else\theshortauthors\fi\fi\fi}}\def\@evenhead{\@oddhead}
\def\@oddfoot{\small\lfoot\ifnum\count0=\startpage\copyright\ \gtp\hfill\else
\gt, Volume \thevolumenumber\ (\thevolumeyear)\hfill\fi}
\def\@evenfoot{\@oddfoot}
\makeatother
\fi

\newwrite\gtoutfile
\long\gdef\makeheadfile{  
{\def\\{, }\def\s{ }
\immediate\openout\gtoutfile head.xxx
\immediate\write\gtoutfile{To: math@arxiv.org}
\immediate\write\gtoutfile{Subject: put or rep NNNNN:pppp}
\immediate\write\gtoutfile{--text follows this line--}
\immediate\write\gtoutfile{Proxy-for: \ifx\theasciiauthors\relax
\theauthors\else\theasciiauthors\fi\s<\ifx\theasciiemail\relax\theemail\else\theasciiemail\fi>}
\immediate\write\gtoutfile{\noexpand\\}
\immediate\write\gtoutfile{Authors: \ifx\theasciiauthors\relax
\theauthors\else\theasciiauthors\fi}
{\def\\{ }\immediate\write\gtoutfile{Title: \ifx\theasciititle\relax
\thetitle\else\theasciititle\fi}}
\immediate\write\gtoutfile{Subj-class: GT or SG or MG etc}
\immediate\write\gtoutfile{MSC-class: \theprimaryclass\ifx\thesecondaryclass\relax\else, \thesecondaryclass\fi}
\immediate\write\gtoutfile{Journal-ref: Geom. Topol. \thevolumenumber
(\thevolumeyear) \startpage-\finishpage}
\immediate\write\gtoutfile{Comments: Published by Geometry and Topology at}
\immediate\write\gtoutfile{\s\s http://www.maths.warwick.ac.uk/gt/GTVol\thevolumenumber/paper\thepapernumber.abs.html}
\immediate\write\gtoutfile{\noexpand\\}
\immediate\write\gtoutfile{}
\ifx\theasciiabstract\relax
\immediate\write\gtoutfile{\theabstract}\else
\immediate\write\gtoutfile{\theasciiabstract}\fi
\immediate\write\gtoutfile{}
\immediate\write\gtoutfile{\noexpand\\}
\immediate\write\gtoutfile{}
\immediate\closeout\gtoutfile}}  

\def\maketitlepage{\maketitlep\makeheadfile}
\let\maketitle\maketitlepage

\lognumber{274}
\volumenumber{7}\papernumber{7}\volumeyear{2003}
\pagenumbers{255}{286}
\received{6 September 2002}
\accepted{12 March 2003}
\published{31 March 2003}
\proposed{Steve Ferry}
\seconded{Benson Farb, Robion Kirby}

\usepackage{amsmath,graphicx,amssymb}
\nocolon

\newtheorem{theorem}{Theorem}[section]

\newtheorem{proposition}[theorem]{Proposition}
\newtheorem{lemma}[theorem]{Lemma}

\theoremstyle{definition}

\newtheorem{example}{Example}

\newtheorem{remark}{Remark}

\numberwithin{equation}{section}

\begin{document}
\title{Manifolds with non-stable fundamental\\groups at infinity, II}
\author{C\thinspace R Guilbault\\F\thinspace C Tinsley}
\coverauthors{C\noexpand\thinspace R Guilbault\\F\noexpand\thinspace C Tinsley}
\asciiauthors{C R Guilbault and F C Tinsley}
\address{Department of Mathematical Sciences, 
University of Wisconsin-Milwaukee\\Milwaukee, Wisconsin 53201, USA}
\email{craigg@uwm.edu, ftinsley@cc.colorado.edu}
\secondaddress{Department of Mathematics, The Colorado College\\Colorado Springs, Colorado 80903, USA}
\asciiaddress{Department of Mathematical Sciences, 
University of Wisconsin-Milwaukee\\Milwaukee, Wisconsin 53201, 
USA\\and\\Department of Mathematics, The Colorado College\\Colorado 
Springs, Colorado 80903, USA}

\begin{abstract}
In this paper we continue an earlier study of ends non-compact manifolds. The
over-arching goal is to investigate and obtain generalizations of Siebenmann's
famous collaring theorem that may be applied to manifolds having non-stable
fundamental group systems at infinity. In this paper we show that, for
manifolds with compact boundary, the condition of inward tameness has
substatial implications for the algebraic topology at infinity. In particular,
every inward tame manifold with compact boundary has stable homology (in all
dimensions) and semistable fundamental group at each of its ends. In contrast,
we also construct examples of this sort which fail to have perfectly
semistable fundamental group at infinity. In doing so, we exhibit the first
known examples of open manifolds that are inward tame and have vanishing Wall
finiteness obstruction at infinity, but are not pseudo-collarable.
\end{abstract}

\maketitlepage

\section{Introduction}

In \cite{Gu1} we presented a program for generalizing Siebenmann's famous
collaring theorem (see \cite{Si}) to include open manifolds with non-stable
fundamental group systems at infinity. To do this, it was first necessary to
generalize the notion of an open collar. Define a manifold $N^{n}$ with
compact boundary to be a \emph{homotopy collar} provided $\partial
N^{n}\hookrightarrow N^{n}$ is a homotopy equivalence. Then define a
\emph{pseudo-collar} to be a homotopy collar which contains arbitrarily small
homotopy collar neighborhoods of infinity. An open manifold is
\emph{pseudo-collarable }if it contains a pseudo-collar neighborhood of
infinity. The main results of our initial investigation may be summarized as follows:

\begin{theorem}
[see \cite{Gu1}]\label{pseudo}Let $M^{n}$ be a one ended $n$-manifold with
compact (possibly empty) boundary. If $M^{n}$ is pseudo-collarable, then

\begin{enumerate}
\item $M^{n}$ is inward tame at infinity,

\item $\pi_{1}(\varepsilon(M^{n}))$ is perfectly semistable, and

\item $\sigma_{\infty}\left(  M^{n}\right)  =0\in\widetilde{K}_{0}\left(
\pi_{1}\left(  \varepsilon(M^{n}\right)  \right)  )$.
\end{enumerate}

\noindent Conversely, for $n\geq7$, if $M^{n}$ satisfies conditions (\rm
(1)--(3) and
$\pi_{2}(\varepsilon\left(  M^{n}\right)  )$ is semistable, then $M^{n}$ is
pseudo-collarable. 
\end{theorem}

\begin{remark}
While it its convenient (and traditional) to focus on one ended manifolds,
this theorem actually applies to all manifolds with compact boundary ---in
particular, to all open manifolds. The key here is that an inward tame
manifold with compact boundary has only finitely many ends---we provide a
proof of this fact in Section \ref{semistability}. Hence, Theorem \ref{pseudo}
may be applied to each end individually. For manifolds with non-compact
boundaries, the situation is quite different. A straight forward
infinite-ended example of this type is given in Section \ref{semistability}. A
more detailed discussion of manifolds with non-compact boundaries will be
provided in \cite{Gu3}.
\end{remark}

The condition of \emph{inward tameness }means (informally) that each
neighborhood of infinity can be pulled into a compact subset of itself. We let
$\pi_{1}(\varepsilon\left(  M^{n}\right)  )$ denote the inverse system of
fundamental groups of neighborhoods of infinity. Such a system is
\emph{semistable} if it is equivalent to a system in which all bonding maps
are surjections. If, in addition, it can be arranged that the kernels of these
bonding maps are perfect groups, then the system is \emph{perfectly
semistable}. The obstruction $\sigma_{\infty}\left(  M^{n}\right)
\in\widetilde{K}_{0}\left(  \pi_{1}\left(  \varepsilon(M^{n}\right)  \right)
)$ vanishes precisely when each (clean) neighborhood of infinity has finite
homotopy type. More precise formulations of these definitions are given in
Section \ref{background}. For a detailed discussion of the structure of
pseudo-collars, along with some useful examples of pseudo-collarable and
non-pseudo-collarable manifolds, the reader is referred to Section 4 of
\cite{Gu1}.

One obvious question suggested by Theorem \ref{pseudo} is whether the
$\pi _{2}$-semistability condition can be omitted from the converse,
ie, whether conditions (1)--(3) are sufficient to guarantee
pseudo-collarability. We are not yet able to resolve that issue. In
this paper, we focus on other questions raised in \cite{Gu1}. The
first asks whether inward tameness implies $\pi_{1}$-semistability;
and the second asks whether inward tameness (possibly combined with
condition 3)) guarantees perfect semistability of $\pi_{1}$.  Thus,
one arrives at the question: Are conditions (1) and (3) sufficient to
ensure pseudo-collarability? Some motivation for this last question is
provided by \cite{CS} where it is shown that these conditions do
indeed characterize pseudo-collarability in Hilbert cube manifolds.

Our first main result provides a positive answer to the $\pi_{1}%
$-semistability question, and more. It shows that---for manifolds with compact
boundary---the inward tameness hypothesis, by itself, has significant
implications for the algebraic topology of that manifold at infinity.

\begin{theorem}
\label{semistable}If an $n$-manifold with compact (possibly empty) boundary is
inward tame at infinity, then it has finitely many ends, each of which has
semistable fundamental group and stable homology in all dimensions.
\end{theorem}

Our second main result provides a negative answer to the pseudo-collarability
question discussed above.

\begin{theorem}
\label{perfect}For $n\geq6$, there exists a one ended open $n$-manifold
$M_{\ast}^{n}$ in which all clean neighborhoods of infinity have finite
homotopy types (hence, $M_{\ast}^{n}$ satisfies conditions (1) and (3) from
above), but which does not have perfectly semistable fundamental group system
at infinity. Thus, $M_{\ast}^{n}$ is not pseudo-collarable.
\end{theorem}

Theorems \ref{semistable} and \ref{perfect} and their proofs are independent.
The first is a very general result that is valid in all dimensions. Its proof
is contained in Section \ref{semistability}. The second involves the
construction of rather specific high-dimensional examples, with a blueprint
being provided by a significant dose of combinatorial group theory. Although
independent, Theorem \ref{semistable} offers crucial guidance on how delicate
such a construction must be. The necessary group theory and the construction
of the examples may be found in Section \ref{construction}. Section
\ref{background} contains the background and definitions needed to read each
of the above. In the final section of this paper we discuss a related open question.

The authors wish to acknowledge Tom Thickstun for some very helpful
discussions.

The first author wishes to acknowledge support from NSF Grant DMS-0072786.

\section{Definitions and Background\label{background}}

This section contains most of the terminology and notation needed in the
remainder of the paper. It is divided into two subsections---the first devoted
to inverse sequences of groups, and the second to the topology of ends of manifolds.

\subsection{Algebra of inverse sequences}

Throughout this section all arrows denote homomorphisms, while arrows of the
type $\twoheadrightarrow$ or $\twoheadleftarrow$ denote surjections. The
symbol $\cong$ denotes isomorphisms.

Let
\[
G_{0}\overset{\lambda_{1}}{\longleftarrow}G_{1}\overset{\lambda_{2}%
}{\longleftarrow}G_{2}\overset{\lambda_{3}}{\longleftarrow}\cdots
\]
be an inverse sequence of groups and homomorphisms. A \emph{subsequence} of
$\left\{  G_{i},\lambda_{i}\right\}  $ is an inverse sequence of the form:
\[
G_{i_{0}}\overset{\lambda_{i_{0}+1}\circ\cdots\circ\lambda_{i_{1}}%
}{\longleftarrow}G_{i_{1}}\overset{\lambda_{i_{1}+1}\circ\cdots\circ
\lambda_{i_{2}}}{\longleftarrow}G_{i_{2}}\overset{\lambda_{i_{2}+1}\circ
\cdots\circ\lambda_{i_{3}}}{\longleftarrow}\cdots
\]
In the future we will denote a composition $\lambda_{i}\circ\cdots\circ
\lambda_{j}$ ($i\leq j$) by $\lambda_{i,j}$.

We say that sequences $\left\{  G_{i},\lambda_{i}\right\}  $ and $\left\{
H_{i},\mu_{i}\right\}  $ are \emph{pro-equivalent} if, after passing to
subsequences, there exists a commuting diagram:
\[%
\begin{array}
[c]{ccccccc}%
G_{i_{0}} & \overset{\lambda_{i_{0}+1,i_{1}}}{\longleftarrow} & G_{i_{1}} &
\overset{\lambda_{i_{1}+1,i_{2}}}{\longleftarrow} & G_{i_{2}} & \overset
{\lambda_{i_{2}+1,i_{3}}}{\longleftarrow} & \cdots\\
& \nwarrow\quad\swarrow &  & \nwarrow\quad\swarrow &  & \nwarrow\quad\swarrow
& \\
& H_{j_{0}} & \overset{\mu_{j_{0}+1,j_{1}}}{\longleftarrow} & H_{j_{1}} &
\overset{\mu_{j_{1}+1,j_{2}}}{\longleftarrow} & H_{j_{2}} & \cdots
\end{array}
\]
Clearly an inverse sequence is pro-equivalent to any of its subsequences. To
avoid tedious notation, we often do not distinguish $\left\{  G_{i}%
,\lambda_{i}\right\}  $ from its subsequences. Instead we simply assume that
$\left\{  G_{i},\lambda_{i}\right\}  $ has the desired properties of a
preferred subsequence---often prefaced by the words ``after passing to a
subsequence and relabelling''.

The \emph{inverse limit }of a sequence $\left\{  G_{i},\lambda_{i}\right\}  $
is a subgroup of $\prod G_{i}$ defined by
\[
\underleftarrow{\lim}\left\{  G_{i},\lambda_{i}\right\}  =\left\{  \left.
\left(  g_{0},g_{1},g_{2},\cdots\right)  \in\prod_{i=0}^{\infty}G_{i}\right|
\lambda_{i}\left(  g_{i}\right)  =g_{i-1}\right\}  .
\]
Notice that for each $i$, there is a \emph{projection homomorphism}
$p_{i}:\underleftarrow{\lim}\left\{  G_{i},\lambda_{i}\right\}  \rightarrow
G_{i}$. It is a standard fact that pro-equivalent inverse sequences have
isomorphic inverse limits.

An inverse sequence $\left\{  G_{i},\lambda_{i}\right\}  $ is \emph{stable} if
it is pro-equivalent to an inverse sequence $\left\{  H_{i},\mu_{i}\right\}  $
for which each $\mu_{i}$ is an isomorphism. Equivalently, $\left\{
G_{i},\lambda_{i}\right\}  $ is stable if, after passing to a subsequence and
relabelling, there is a commutative diagram of the form
\begin{equation}%
\begin{array}
[c]{ccccccccc}%
G_{0} & \overset{\lambda_{1}}{\longleftarrow} & G_{1} & \overset{\lambda_{2}%
}{\longleftarrow} & G_{2} & \overset{\lambda_{3}}{\longleftarrow} & G_{3} &
\overset{\lambda_{4}}{\longleftarrow} & \cdots\\
& \nwarrow\quad\swarrow &  & \nwarrow\quad\swarrow &  & \nwarrow\quad\swarrow
&  &  & \\
& im(\lambda_{1}) & \longleftarrow & im(\lambda_{2}) & \longleftarrow &
im(\lambda_{3}) & \longleftarrow & \cdots &
\end{array}
\tag{$\ast$}%
\end{equation}
where each bonding map in the bottom row (obtained by restricting the
corresponding $\lambda_{i}$) is an isomorphism. If $\left\{  H_{i},\mu
_{i}\right\}  $ can be chosen so that each $\mu_{i}$ is an epimorphism, we say
that our inverse sequence is \emph{semistable }(or \emph{Mittag-Leffler,
}or\emph{\ pro-epimorphic}). In this case, it can be arranged that the
restriction maps in the bottom row of ($\ast$) are epimorphisms. Similarly, if
$\left\{  H_{i},\mu_{i}\right\}  $ can be chosen so that each $\mu_{i}$ is a
monomorphism, we say that our inverse sequence is \emph{pro-monomorphic}; it
can then be arranged that the restriction maps in the bottom row of ($\ast$)
are monomorphisms. It is easy to see that an inverse sequence that is
semistable and pro-monomorphic is stable.

Recall that a \emph{commutator} element of a group $H$ is an element of the
form $x^{-1}y^{-1}xy$ where $x,y\in H$; and the \emph{commutator subgroup} of
$H,$ denoted $\left[  H,H\right]  $, is the subgroup generated by all of its
commutators. The group $H$ is \emph{perfect} if $\left[  H,H\right]  =H$. An
inverse sequence of groups is \emph{perfectly semistable} if it is
pro-equivalent to an inverse sequence
\[
G_{0}\overset{\lambda_{1}}{\twoheadleftarrow}G_{1}\overset{\lambda_{2}%
}{\twoheadleftarrow}G_{2}\overset{\lambda_{3}}{\twoheadleftarrow}\cdots
\]
of finitely presentable groups and surjections where each $\ker\left(
\lambda_{i}\right)  $ perfect. The following shows that inverse sequences of
this type behave well under passage to subsequences.

\begin{lemma}
\label{L1}A composition of surjective group homomorphisms, each having perfect
kernels, has perfect kernel. Thus, if an inverse sequence of surjective group
homomorphisms has the property that the kernel of each bonding map is perfect,
then each of its subsequences also has this property.
\end{lemma}

\begin{proof}
See Lemma 1 of \cite{Gu1}.
\end{proof}

For later use, we record an easy but crucial property of perfect groups.

\begin{lemma}
\label{L2}If $f:G\twoheadrightarrow H$ is a surjective group homomorphism and
$G$ is perfect, then $H$ is perfect.
\end{lemma}

\begin{proof}
The image of each commutator from $G$ is a commutator in $H.$
\end{proof}

We conclude this section with a technical result that will be needed later.
Compare to the well-known Five Lemma from homological algebra.

\begin{lemma}
\label{dydak}Assume the following commutative diagram of five inverse
sequences:%
\[%
\begin{array}
[c]{ccccccccc}%
\vdots &  & \vdots &  & \vdots &  & \vdots &  & \vdots\\
\downarrow &  & \downarrow &  & \downarrow &  & \downarrow &  & \downarrow\\
A_{2} & \rightarrow & B_{2} & \rightarrow & C_{2} & \rightarrow & D_{2} &
\rightarrow & E_{2}\\
\downarrow &  & \downarrow &  & \downarrow &  & \downarrow &  & \downarrow\\
A_{1} & \rightarrow & B_{1} & \rightarrow & C_{1} & \rightarrow & D_{1} &
\rightarrow & E_{1}\\
\downarrow &  & \downarrow &  & \downarrow &  & \downarrow &  & \downarrow\\
A_{0} & \rightarrow & B_{0} & \rightarrow & C_{0} & \rightarrow & D_{0} &
\rightarrow & E_{0}%
\end{array}
\]
If each row is exact and the inverse sequences $\left\{  A_{i}\right\}  $,
$\left\{  B_{i}\right\}  $, $\left\{  D_{i}\right\}  $, and $\left\{
E_{i}\right\}  $ are stable, then so is $\left\{  C_{i}\right\}  $.

\begin{proof}
The proof is by an elementary but intricate diagram chase. See Lemmas 2.1 and
2.2 of \cite{Dy}.
\end{proof}
\end{lemma}

\subsection{Topology of ends of manifolds\label{topology}}

In this paper, the term \emph{manifold} means \emph{manifold with (possibly
empty) boundary}. A manifold is \emph{open} if it is non-compact and has no
boundary. For convenience, all manifolds are assumed to be PL. Analogous
results may be obtained for smooth or topological manifolds in the usual ways.

Let $M^{n}$ be a manifold with compact (possibly empty) boundary. A set
$N\subset M^{n}$ is a \emph{neighborhood of infinity} if $\overline{M^{n}-N}$
is compact. A neighborhood of infinity $N$ is \emph{clean} if

\begin{itemize}
\item $N$ is a closed subset of $M^{n}$,

\item $N\cap\partial M^{n}=\emptyset$, and

\item $N$ is a codimension 0 submanifold of $M^{n}$ with bicollared boundary.
\end{itemize}

\noindent It is easy to see that each neighborhood of infinity contains a
clean neighborhood of infinity.

\begin{remark}
We have taken advantage of the compact boundary by requiring that clean
neighborhoods of infinity be disjoint from $\partial M^{n}$. In the case of
non-compact boundary, a slightly more delicate definition is required.
\end{remark}

We say that $M^{n}$ \emph{has }$k$ \emph{ends }if it contains a compactum $C$
such that, for every compactum $D$ with $C\subset D$, $M^{n}-D$ has exactly
$k$ unbounded components, ie, $k$ components with noncompact closures. When
$k$ exists, it is uniquely determined; if $k$ does not exist, we say $M^{n}$
has \emph{infinitely many ends}.

If $M^{n}$ has compact boundary and is $k$-ended, then $M^{n}$ contains a
clean neighborhood of infinity $N$ that consists of $k$ connected components,
each of which is a one ended manifold with compact boundary. Therefore, when
studying manifolds (or other spaces) having finitely many ends, it suffices to
understand the \emph{one ended} situation. In this paper, we are primarily
concerned with manifolds possessing finitely many ends (See Theorem
\ref{semistable} or Prop. \ref{1}), and thus, we frequently restrict our
attention to the one ended case.

A connected clean neighborhood of infinity with connected boundary is called a
\emph{0-neighborhood of infinity}. If $N$ is clean and connected but has more
than one boundary component, we may choose a finite collection of disjoint
properly embedded arcs in $N$ that connect these components. Deleting from $N$
the interiors of regular neighborhoods of these arcs produces a $0$%
-neighborhood of infinity $N_{0}\subset N$.

A nested sequence $N_{0}\supset N_{1}\supset N_{2}\supset\cdots$ of
neighborhoods of infinity is \emph{cofinal }if $\bigcap_{i=0}^{\infty}%
N_{i}=\emptyset$. For any one ended manifold $M^{n}$ with compact boundary,
one may easily obtain a cofinal sequence of $0$-neighborhoods of infinity.

We say that $M^{n\text{ }}$is \emph{inward tame }at infinity if, for
arbitrarily small neighborhoods of infinity $N$, there exist homotopies
$H:N\times\left[  0,1\right]  \rightarrow N$ such that $H_{0}=id_{N}$ and
$\overline{H_{1}\left(  N\right)  }$ is compact. Thus inward tameness means
each neighborhood of infinity can be pulled into a compact subset of itself.
In this situation, the $H$'s will be referred to as \emph{taming homotopies.}

Recall that a complex $X$ is \emph{finitely dominated} if there exists a
finite complex $K$ and maps $u:X\rightarrow K$ and $d:K\rightarrow X$ such
that $d\circ u\simeq id_{X}$. The following lemma uses this notion to offer
equivalent formulations of ``inward tameness''.

\begin{lemma}
For a manifold $M^{n}$, the following are equivalent.

\begin{enumerate}
\item $M^{n}$ is inward tame at infinity.

\item Each clean neighborhood of infinity in $M^{n}$ is finitely dominated.

\item For each cofinal sequence $\left\{  N_{i}\right\}  $ of clean
neighborhoods of infinity, the inverse sequence
\[
N_{0}\overset{j_{1}}{\hookleftarrow}N_{1}\overset{j_{2}}{\hookleftarrow}%
N_{2}\overset{j_{3}}{\hookleftarrow}\cdots
\]
is pro-homotopy equivalent to an inverse sequence of finite polyhedra.
\end{enumerate}

\begin{proof}
To see that (1) implies (2), let $N$ be a clean neighborhood of infinity and
$H:N\times\left[  0,1\right]  \rightarrow N$ a taming homotopy. Let $K$ be a
polyhedral subset of $N$ that contains $\overline{H_{1}\left(  N\right)  }$.
If $u:N\rightarrow K$ is obtained by restricting the range of $H_{1}$ and
$d:K\hookrightarrow N$, then $d\circ u=H_{1}\simeq id_{N}$, so $N$ is finitely dominated.

To see that 2) implies 3), choose for each $N_{i}$ a finite polyhedron $K_{i}$
and maps $u_{i}:N_{i}\rightarrow K_{i}$ and $d_{i}:K_{i}\rightarrow N_{i}$
such that $d_{i}\circ u_{i}\simeq id_{N_{i}}$. For each $i\geq1$, let
$f_{i}=u_{i-1}\circ j_{i}$ and $g_{i}=f_{i}\circ d_{i}$. Since $d_{i-1}\circ
f_{i}=d_{i-1}\circ u_{i-1}\circ j_{i}\simeq id_{N_{i-1}}\circ j_{i}=j_{i}$,
the diagram%
\[%
\begin{array}
[c]{ccccccccc}%
N_{0} & \overset{j_{1}}{\hookleftarrow} & N_{1} & \overset{j_{2}%
}{\hookleftarrow} & N_{2} & \overset{j_{3}}{\hookleftarrow} & N_{3} &
\overset{j_{4}}{\hookleftarrow} & \cdots\\
& \overset{\quad d_{0}}{\nwarrow}\quad\overset{f_{1\quad}}{\swarrow} &  &
\overset{\quad d_{1}}{\nwarrow}\quad\overset{f_{2\quad}}{\swarrow} &  &
\overset{\quad d_{2}}{\nwarrow}\quad\overset{f_{3\quad}}{\swarrow} &  &  & \\
& K_{0} & \overset{g_{1}}{\longleftarrow} & K_{1} & \overset{g_{2}%
}{\longleftarrow} & K_{2} & \overset{g_{3}}{\longleftarrow} & \cdots &
\end{array}
\]
commutes up to homotopy, so (by definition) the two inverse sequences are
pro-homotopy equivalent.

Lastly, we assume the existence of a homotopy commutative diagram as pictured
above for some cofinal sequence of clean neighborhoods of infinity and some
inverse sequence of finite polyhedra. We show that for each $i\geq1$, there is
a taming homotopy for $N_{i}$. By hypothesis, $d_{i}\circ f_{i+1}\simeq
j_{i+1}$. Extend $j_{i+1}$ to $id_{N_{i}}$, then apply the homotopy extension
property (see \cite[pp.14-15]{Ha}) for the pair $\left(  N_{i},N_{i+1}\right)
$ to obtain $H:N_{i}\times\left[  0,1\right]  \rightarrow N_{i}$ with
$H_{0}=id_{N_{i}}$ and $\left.  H_{1}\right|  _{N_{i+1}}=d_{i}\circ f_{i+1}$.
Now,
\[
H_{1}\left(  N_{i}\right)  =H_{1}\left(  N_{i}-N_{i+1}\right)  \cup
H_{1}\left(  N_{i+1}\right)  \subset H_{1}\left(  \overline{N_{i}-N_{i+1}%
}\right)  \cup d_{i}\left(  K_{i}\right)  ,
\]
so $\overline{H_{1}\left(  N_{i}\right)  }$ is compact, and $H$ is the desired
taming homotopy.
\end{proof}
\end{lemma}

Given a nested cofinal sequence $\left\{  N_{i}\right\}  _{i=0}^{\infty}$ of
connected neighborhoods of infinity, base points $p_{i}\in N_{i}$, and paths
$\alpha_{i}\subset N_{i}$ connecting $p_{i}$ to $p_{i+1}$, we obtain an
inverse sequence:
\[
\pi_{1}\left(  N_{0},p_{0}\right)  \overset{\lambda_{1}}{\longleftarrow}%
\pi_{1}\left(  N_{1},p_{1}\right)  \overset{\lambda_{2}}{\longleftarrow}%
\pi_{1}\left(  N_{2},p_{2}\right)  \overset{\lambda_{3}}{\longleftarrow}%
\cdots
\]
Here, each $\lambda_{i+1}:\pi_{1}\left(  N_{i+1},p_{i+1}\right)
\rightarrow\pi_{1}\left(  N_{i},p_{i}\right)  $ is the homomorphism induced by
inclusion followed by the change of base point isomorphism determined by
$\alpha_{i}$. The obvious singular ray obtained by piecing together the
$\alpha_{i}$'s is often referred to as the \emph{base ray }for the inverse
sequence. Provided the sequence is semistable, one can show that its
pro-equivalence class does not depend on any of the choices made above. We
refer to the pro-equivalence class of this sequence as the \emph{fundamental
group system at infinity} for $M^{n}$ and denote it by $\pi_{1}\left(
\varepsilon\left(  M^{n}\right)  \right)  $. (In the absence of semistability,
the pro-equivalence class of the inverse sequence depends on the choice of
base ray, and hence, this choice becomes part of the data.) It is easy to see
how the same procedure may also be used to define $\pi_{k}\left(
\varepsilon\left(  M^{n}\right)  \right)  $ for $k>1$.

For any coefficient ring $R$ and any integer $j\geq0$, a similar procedure
yields an inverse sequence
\[
H_{j}\left(  N_{0};R\right)  \overset{\lambda_{1}}{\longleftarrow}H_{j}\left(
N_{1};R\right)  \overset{\lambda_{2}}{\longleftarrow}H_{j}\left(
N_{2};R\right)  \overset{\lambda_{3}}{\longleftarrow}\cdots
\]
where each $\lambda_{i}$ is induced by inclusion---here, no base points or
rays are needed. We refer to the pro-equivalence class of this sequence as the
\emph{j}$^{th}$\emph{\ homology at infinity }for $M^{n}$ with $R$-coefficients
and denote it by $H_{j}\left(  \varepsilon\left(  M^{n}\right)  ;R\right)  $.

In \cite{Wll}, Wall shows that each finitely dominated connected space $X$
determines a well-defined element $\sigma\left(  X\right)  $ lying in
$\widetilde{K}_{0}\left(  \mathbb{Z}\left[  \pi_{1}X\right]  \right)  $ (the
group of stable equivalence classes of finitely generated projective
$\mathbb{Z}\left[  \pi_{1}X\right]  $-modules under the operation induced by
direct sum) that vanishes if and only if $X$ has the homotopy type of a finite
complex. Given a nested cofinal sequence $\left\{  N_{i}\right\}
_{i=0}^{\infty}$ of connected clean neighborhoods of infinity in an inward
tame manifold $M^{n}$, we have a Wall obstruction $\sigma(N_{i})$ for each
$i$. These may be combined into a single obstruction
\[
\begin{split}
\sigma_{\infty}(M^{n})&=\left(  -1\right)  ^{n}(\sigma(N_{0}),\sigma
(N_{1}),\sigma(N_{2}),\cdots)\\
&\in\widetilde{K}_{0}\left(  \pi_{1}\left(
\varepsilon\left(  M^{n}\right)  \right)  \right)  \equiv\underleftarrow{\lim
}\widetilde{K}_{0}\left(  \mathbb{Z}\left[  \pi_{1}N_{i}\right]  \right)
\end{split}
\]
that is well-defined and which vanishes if and only if each clean neighborhood
of infinity in $M^{n}$ has finite homotopy type. See \cite{CS} for details.

We close this section with a known result from the topology of manifolds. Its
proof is short and its importance is easily seen when one considers the
``one-sided $h$-cobordism'' $\left(  W,\partial N,\partial N^{\prime}\right)
$ that occurs naturally when $N^{\prime}$ is a homotopy collar contained in
the interior of another homotopy collar $N$ and $W=\overline{N-N^{\prime}}$.
In particular, this result explains why pseudo-collarable manifolds must have
perfectly semistable fundamental groups at their ends. Additional details may
be found in Section 4 of \cite{Gu1}.

\begin{theorem}
Let $\left(  W^{n},P,Q\right)  $ be a compact connected cobordism between
closed $\left(  n-1\right)  $-manifolds with the property that
$P\hookrightarrow W^{n}$ is a homotopy equivalence. Then the inclusion induced
map $i_{\#}:\pi_{1}(Q)\rightarrow\pi_{1}(W^{n})$ is surjective and has perfect kernel.

\begin{proof}
Let $p:\widetilde{W}\rightarrow W^{n}$ be the universal covering projection,
$\widetilde{P}=p^{-1}\left(  P\right)  $, and $\widehat{Q}=p^{-1}\left(
Q\right)  $. By Poincar\'{e} duality for non-compact manifolds, $$H_{k}\left(
\widetilde{W},\widehat{Q};\mathbb{Z}\right)  \cong H_{c}^{n-k}\left(
\widetilde{W},\widetilde{P};\mathbb{Z}\right),$$ where cohomology is with
compact supports. Since $\widetilde{P}\hookrightarrow\widetilde{W}$ is a
proper homotopy equivalence, all of these relative cohomology groups vanish.
It follows that $H_{1}\left(  \widetilde{W},\widehat{Q};\mathbb{Z}\right)
=0$, so by the long exact sequence for $\left(  \widetilde{W},\widehat
{Q}\right)  $, $\widetilde{H}_{0}\left(  \widehat{Q};\mathbb{Z}\right)  =0$;
therefore $\widehat{Q}$ is connected. By covering space theory, the components
of $\widehat{Q}$ are in one-to-one correspondence with the cosets of
$i_{\#}\left(  \pi_{1}(Q)\right)  $ in $\pi_{1}(W^{n})$, so $i_{\#}$ is
surjective. Similarly, $H_{2}\left(  \widetilde{W},\widehat{Q};\mathbb{Z}%
\right)  =0$, and since $\widetilde{W}$ is simply connected, the long exact
sequence for $\left(  \widetilde{W},\widehat{Q}\right)  $ shows that
$H_{1}\left(  \widehat{Q};\mathbb{Z}\right)  =0$. This implies that $\pi
_{1}\left(  \widehat{Q}\right)  $ is a perfect group, and covering space
theory tell us that $\pi_{1}\left(  \widehat{Q}\right)  \cong\ker\left(
i_{\#}\right)  $.
\end{proof}
\end{theorem}

\section{Inward tameness, $\pi_{1}$-semistability, and $H_{\ast}%
$-stability\label{semistability}}

The theme of this section is that---for manifolds with compact (possibly
empty) boundary---inward tameness, by itself, has some significant
consequences. In particular, an inward tame manifold of this type has:

\begin{itemize}
\item finitely many ends,

\item semistable fundamental group at each of these ends, and

\item stable (finitely generated) homology at infinity in all dimensions.
\end{itemize}

\noindent The first of these properties is known; for completeness, we will
provide a proof. The second property answers a question posed in \cite{Gu1}. A
stronger conclusion of $\pi_{1}$-stability is not possible, as can be seen in
the exotic universal covering spaces constructed in \cite{Da}. (See Example 3
of \cite{Gu1} for a discussion.) Somewhat surprisingly, inward tameness
\emph{does} imply stability at infinity for homology in the situation at hand.

It is worth noting that, under slightly weaker hypotheses, none of these
properties holds. We provide some simple examples of locally finite complexes,
and polyhedral manifolds (with non-compact boundaries) that violate each of
the above.

\begin{example}
Let $E$ denote a wedge of two circles. Then the universal cover $\widetilde
{E}$ of $E$ is an inward tame $1$-complex with infinitely many ends.
\end{example}

\begin{example}
Let $f:(S^{1},\ast)\rightarrow(S^{1},\ast)$ be degree $2$ map, and let $X$ be
the ``inverse mapping telescope'' of the system:%
\[
S^{1}\overset{f}{\longleftarrow}S^{1}\overset{f}{\longleftarrow}S^{1}%
\overset{f}{\longleftarrow}\cdots
\]
Assemble a base ray from the mapping cylinder arcs corresponding to the base
point $\ast$. It is easy to see that $X$ is inward tame and that $\pi
_{1}\left(  \varepsilon(X\right)  )$ is represented by the system%
\[
\mathbb{Z}\overset{\times2}{\longleftarrow}\mathbb{Z}\overset{\times
2}{\longleftarrow}\mathbb{Z}\overset{\times2}{\longleftarrow}\cdots%
\]
which is not semistable. Hence, $\pi_{1}$-semistability does not follow from
inward tameness for one ended complexes. This example also shows that inward
tame complexes needn't have stable $H_{1}\left(  \varepsilon(X\right)
;\mathbb{Z})$.
\end{example}

\begin{example}
More generally, if
\[
K_{0}\overset{f_{1}}{\longleftarrow}K_{1}\overset{f_{2}}{\longleftarrow}%
K_{3}\overset{f_{3}}{\longleftarrow}\cdots
\]
is an inverse sequence of finite polyhedra, then the inverse mapping telescope
$Y$ of this sequence is inward tame. By choosing the polyhedra and the bonding
maps appropriately, we can obtain virtually any desired behavior in $\pi
_{1}\left(  \varepsilon(Y\right)  )$ and $H_{k}\left(  \varepsilon(Y\right)
;\mathbb{Z}).$
\end{example}

\begin{example}
By properly embedding the above complexes in $\mathbb{R}^{n}$ and letting
$M^{n}$ be a regular neighborhood, we may obtain inward tame manifold examples
with similar bad behavior at infinity. Of course, $M^{n}$ will have noncompact boundary.
\end{example}

We are now ready to prove Theorem \ref{semistable}. This will be done with a
sequence of three propositions---one for each of the bulleted items listed
above. The first is the simplest and may be deduced from Theorem 1.10 of
\cite{Si}. It could also be obtained later, as a corollary of Proposition
\ref{3}. However, Proposition \ref{3} and its proof become cleaner if we
obtain this result first. The proof is short and rather intuitive.

\begin{proposition}
\label{1}Let $M^{n}$ be an $n$-manifold with compact boundary that is inward
tame at infinity. Then $M^{n}$ has finitely many ends. More specifically, the
number of ends is less than or equal to $rank\left(  H_{n-1}\left(
M^{n};\mathbb{Z}_{2}\right)  \right)  +1$. (See the remark below.)
\end{proposition}

\begin{proof}
Inward tameness implies that each clean neighborhood of infinity (including
$M^{n}$ itself) is finitely dominated and hence, has finitely generated
homology in all dimensions. We'll show that $M^{n}$ has at most $k_{0}+1$
ends, where $k_{0}=rank\left(  H_{n-1}\left(  M^{n};\mathbb{Z}_{2}\right)
\right)  .$

Let $N$ be an clean neighborhood of infinity, each of whose components is
noncompact. Since $H_{0}\left(  N;\mathbb{Z}_{2}\right)  $ has finite rank,
there are finitely many of these components $\left\{  N_{i}\right\}
_{i=1}^{p} $. Our theorem follows if we can show that $p$ is bounded by
$k_{0}+1$.

Using techniques described in Section \ref{topology}, we may assume that
$\partial N_{i}$ is non-empty and connected for all $i$. Then, from the long
exact sequence for the pair $\left(  N_{i},\partial N_{i}\right)  $, we may
deduce that for each $i$, $rank\left(  H_{n-1}\left(  N_{i};\mathbb{Z}%
_{2}\right)  \right)  \geq1$. Hence, $rank\left(  H_{n-1}\left(
N;\mathbb{Z}_{2}\right)  \right)  \geq p$

Let $C=\overline{M^{n}-N}$. Then $C$ is a compact codimension $0$ submanifold
of $M^{n}$, and its boundary consists of the disjoint union of $\partial
M^{n}$ with $\partial N.$ Thus, $rank\left(  H_{n-1}\left(  \partial
C;\mathbb{Z}_{2}\right)  \right)  =p+q$, where $q$ is the number of components
in $\partial M^{n}$. From the long exact sequence for the pair $\left(
C,\partial C\right)  $ we may conclude that $rank\left(  H_{n-1}\left(
C;\mathbb{Z}_{2}\right)  \right)  \geq p+q-1$.

Now consider the following Mayer-Vietoris sequence:%
\[%
\begin{array}
[c]{ccccc}%
\rightarrow H_{n-1}\left(  \partial N;\mathbb{Z}_{2}\right)  & \rightarrow &
H_{n-1}\left(  C;\mathbb{Z}_{2}\right)  \oplus H_{n-1}\left(  N;\mathbb{Z}%
_{2}\right)  & \rightarrow & H_{n-1}\left(  M^{n};\mathbb{Z}_{2}\right)  
\rightarrow\\
 \shortparallel &  &  &  & \shortparallel \\
 \bigoplus_{i=1}^{p}\mathbb{Z}_{2} &  &  &  & \bigoplus_{i=1}^{k_{0}%
}\mathbb{Z}_{2} 
\end{array}
\]

Since $\mathbb{Z}_{2}$ is a field, exactness implies that the rank of the
middle term is no greater than the sum of the ranks of the first and third
terms. The first summand of the middle term has rank $\geq p+q-1$ and the
second summand has rank $\geq p$. Hence $2p+q-1\leq p+k_{0}$. It follows that
$p\leq k_{0}+1$.
\end{proof}

\begin{remark}
The number of ends of $M^{n}$ may be less than $rank\left(  H_{n-1}\left(
M^{n};\mathbb{Z}_{2}\right)  \right)$  $+1$. Indeed, by ``connect summing''
copies of $S^{n-1}\times S^{1}$ to $\mathbb{R}^{n}$, one can make the
difference between these numbers arbitrarily large. The issue is that some
generators of $H_{n-1}\left(  M^{n};\mathbb{Z}_{2}\right)  $ do not ``split
off an end''. To obtain strict equality one should add $1$ to the rank of the
kernel of
\[
\lambda:H_{n-1}\left(  M^{n};\mathbb{Z}_{2}\right)  \rightarrow H_{n-1}%
^{lf}\left(  M^{n};\mathbb{Z}_{2}\right)
\]
where $H^{lf}$ denotes homology based on locally finite chains.
\end{remark}

Before proving the remaining two propositions, we fix some notation and
describe a ``homotopy refinement procedure'' that will be applied in each of
the proofs. As noted earlier, (by applying Proposition \ref{1}) it suffices to
consider the one ended case, so for the remainder of this section, $M^{n}$ is
a one ended inward tame manifold with compact boundary.

Let $\left\{  N_{i}\right\}  _{i=0}^{\infty}$ be a nested cofinal sequence of
$0$-neighborhoods of infinity and, for each $i\geq0$, let $A_{i}%
=N_{i}-int(N_{i+1})$. By inward tameness, we may (after passing to a
subsequence and relabelling) assume that (for each $i\geq0$) there exists a
taming homotopy $H^{i}:N_{i}\times\left[  0,1\right]  \rightarrow N_{i}$ satisfying:

\begin{enumerate}
\item[i)] $H_{0}^{i}=id_{N_{i}}$,

\item[ii)] $H^{i}$ is fixed on $\partial N_{i}$, and

\item[iii)] $H_{1}^{i}\left(  N_{i}\right)  \subset A_{i}-\partial N_{i+1}$.
\end{enumerate}

Choose a proper embedding $r:[0,\infty)\rightarrow N_{0}$ so that, for each
$i$, $r\left(  [i,\infty)\right)  \subset N_{i}$ and so that the image ray
$R_{0}$ intersects each $\partial N_{i}$ transversely at the single point
$p_{i}=r\left(  i\right)  $. For $i\geq0$, let $R_{i}=r\left(  [i,\infty
)\right)  \subset N_{i}$; and let $\alpha_{i}$ denote the arc $r\left(
\left[  i,i+1\right]  \right)  $ in $A_{i}$ from $p_{i}$ to $p_{i+1}$. In
addition, choose an embedding $t:B^{n-1}\times\lbrack0,\infty)\rightarrow
N_{0}$ whose image $T_{0}$ is a regular neighborhood of $R_{0}$, such that
$\left.  t\right|  _{\left\{  \overline{0}\right\}  \times\lbrack0,\infty)}%
=r$, and so that, for each $i$, $T_{0}$ intersects $\partial N_{i}$ precisely
in the $\left(  n-1\right)  $-disk $D_{i}=t(B^{n-1}\times\left\{  i\right\}
)$. Let $B^{\prime}\subset int\left(  B^{n-1}\right)  $ be an $\left(
n-1\right)  $-ball containing $\overline{0}$, $T_{0}^{\prime}=t\left(
B^{\prime}\times\lbrack0,\infty)\right)  $ and $D_{i}^{\prime}=t(B^{\prime
}\times\left\{  i\right\}  )$. 
\begin{figure}
[ht!]
\begin{center}
\includegraphics[
width=.8\hsize
]%
{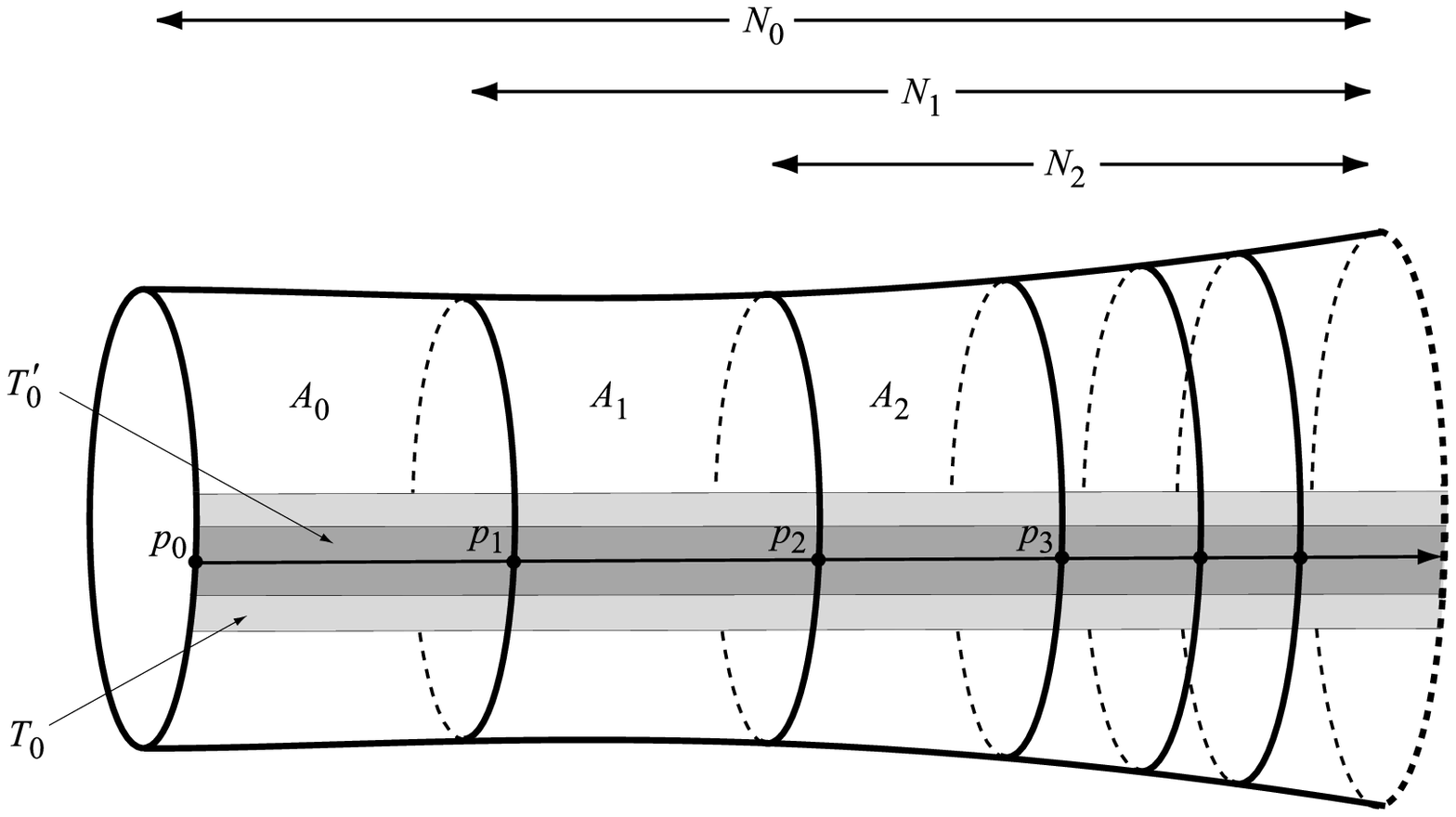}%
\caption{ }%
\end{center}
\end{figure}%
Then, for each $i\geq0$, $T_{i}=t(B^{n-1}%
\times\lbrack i,\infty))$ and $T_{i}^{\prime}=t\left(  B^{\prime}\times\lbrack
i,\infty)\right)  $ are regular neighborhoods of $R_{i}$ in $N_{i}$
intersecting $\partial N_{i}$ in $D_{i}$ and $D_{i}^{\prime}$, respectively.
See Figure 1.%

We now show how to refine each $H^{i}$ so that it respects the ``base ray''
$R_{i}$ and acts in a particularly nice manner on and over $T_{i}^{\prime}$.
Let $j^{i}:(B^{n-1}\times\lbrack i,\infty))\times\left[  0,1\right]
\rightarrow B^{n-1}\times\lbrack i,\infty)$ be a strong deformation retraction
onto $\partial\left(  B^{n-1}\times\lbrack i,\infty)\right)  $ with the
following properties:

\begin{enumerate}
\item[a)] On $B^{\prime}\times\lbrack i,\infty)$, $j^{i}$ is the ``radial''
deformation retraction onto $B^{\prime}\times\left\{  i\right\}  $ given by
$(\left(  b,s\right)  ,u)\mapsto\left(  b,s+u\left(  i-s\right)  \right)  $.

\item[b)] For $(b,s)\notin B^{\prime}\times\lbrack i,\infty)$, the track
$j^{i}((b,s)\times\left[  0,1\right]  )$ of $\left(  b,s\right)  $ does not
intersect $B^{\prime}\times\lbrack i,\infty)$.

\item[c)] The radial component of each track of $j^{i}$ is non-increasing,
ie, if $u_{1}\leq u_{2}$ then $p(j^{i}\left(  b,s,u_{2}\right)  )\leq
p\left(  j^{i}\left(  b,s,u_{1}\right)  \right)  $ where $p$ is projection
onto $[i,\infty)$.
\end{enumerate}

Figure 2 represents $j^{i}$, wherein tracks of $j^{i}$ are meant to
follow the indicated flow lines.

\begin{figure}
[ht!]
\begin{center}
\includegraphics[
width=.9\hsize
]%
{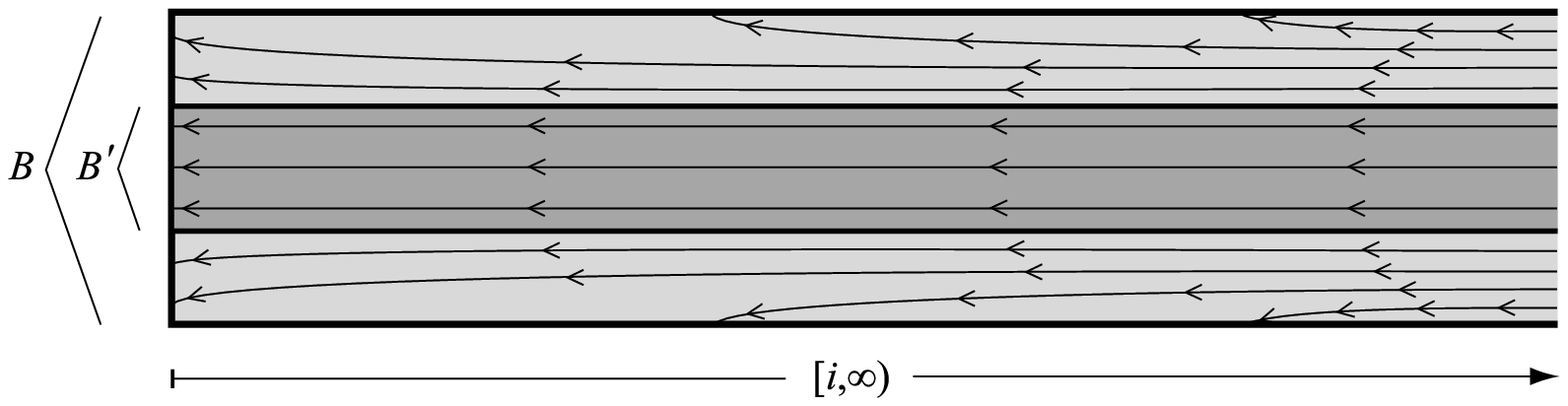}%
\caption{ }%
\end{center}
\end{figure}

Define $J^{i}:N_{i}\times\left[  0,1\right]  \rightarrow N_{i}$ to be $t\circ
j^{i}\circ(t^{-1}\times id)$ on $T_{i}$ and the identity outside of $T_{i}$.
Then $J^{i}$ is a strong deformation retraction of $N_{i}$ onto $N_{i}%
-t\left(  \mathring{B}^{n-1}\times\left(  i,\infty\right)  \right)  $. Define
$K^{i}:N_{i}\times\left[  0,1\right]  \rightarrow N_{i}$ as follows:%
\[
K^{i}\left(  x,t\right)  =\left\{
\begin{array}
[c]{cc}%
J^{i}(x,2t) & 0\leq t\leq\frac{1}{2}\\
J_{1}^{i}(H^{i}\left(  J^{i}\left(  x,1\right)  ,2t-1\right)  ) & \frac{1}%
{2}\leq t\leq1
\end{array}
\right.  \text{.}%
\]
This homotopy retains the obvious analogs of properties i)-iii). In addition,
we have

\begin{enumerate}
\item[iv)] $K^{i}$ acts in a canonical manner on $T_{i}^{\prime}$, and

\item[v)] tracks of points outside of $T_{i}^{\prime}$ do not pass through the
interior of $T_{i}^{\prime}$.
\end{enumerate}

\begin{proposition}
\label{2}Every one ended inward tame $n$-manifold $M^{n}$ with compact
boundary has semistable fundamental group at infinity.
\end{proposition}

\begin{proof}
For convenience, assume that $n\geq3$. For $n=2$ the result may be obtained by
applying well-known structure theorems for $2$-manifolds, or by modifying our
proof slightly.

Let $\left\{  N_{i}\right\}  _{i=0}^{\infty}$ be a nested cofinal sequence of
$0$-neighborhoods of infinity with refined taming homotopies $\left\{
K^{i}\right\}  _{i=0}^{\infty}$ as constructed above. Other choices and labels
are also carried over from above. Note that, for each $i$,
\begin{equation}
\pi_{1}\left(  A_{i},p_{i}\right)  \rightarrow\pi_{1}\left(  N_{i}%
,p_{i}\right)  \text{ is surjective} \tag{$\dagger$}%
\end{equation}

We will show that (for each $i\geq0$) each loop in $N_{i+1}$ based at
$p_{i+1}$ can be pushed (rel $\alpha_{i+1}$) to a loop in $N_{i+2}$ based at
$p_{i+2}$ via a homotopy contained in $N_{i}$. This implies the existence of a
diagram of type ($\ast$) from section \ref{background} for which the bonding
homomorphisms in the bottom row are surjective---and thus, $\pi_{1}$-semistability.

(\textbf{Note}\qua
In performing this push, the ``rel $\alpha_{i+1}$'' requirement
is crucial. The ability to push loops from $N_{i+1}$ into $N_{i+2}$ via a
homotopy contained in $N_{i}$---without regards to basepoints---would yield
another well-known, but strictly weaker, property called \emph{end
1-movability}. See \cite{BT} for a discussion. Much of the homotopy refinement
process described earlier is aimed at obtaining control over the tracks of the
base points.)

Let $\tau$ be a chosen loop in $N_{i+1}$ based at $p_{i+1}$. By ($\dagger$),
we may assume that $\tau\subset A_{i+1}-\partial N_{i+2}$. Let $L^{i}:\partial
N_{i+2}\times\left[  0,1\right]  \rightarrow N_{i}$ be the restriction of
$K^{i}$. Note that $L^{i}\left(  \partial N_{i+2}\times\left\{  1\right\}
\right)  \subset A_{i}-\partial N_{i+1}$ and that, by condition (iv) above,
$L^{i}$ takes $D_{i+2}^{\prime}\times\left[  0,\frac{1}{2}\right]  $
homeomorphically onto $\overline{T_{i}^{\prime}-T_{i+2}^{\prime}}$ with
$D_{i+2}^{\prime}\times\left[  \frac{1}{2},1\right]  $ being flattened onto
$D_{i}$. In addition, $L^{i}(\left\{  p_{i+2}\right\}  \times\left[
0,\frac{1}{4}\right]  )=\alpha_{i+1}$, $L^{i}(\left\{  p_{i+2}\right\}
\times\left[  \frac{1}{4},\frac{1}{2}\right]  )=\alpha_{i}$, and
$L^{i}(p_{i+2},\frac{1}{4})=p_{i+1}$. Without changing its values on
$(\partial N_{i+2}\times\left\{  0\right\}  )\cup(D_{i+2}^{\prime}%
\times\left[  0,\frac{1}{2}\right]  )$, we may adjust $L^{i} $ so that it is a
non-degenerate PL mapping. In particular, we may choose triangulations
$\Gamma_{1}$ and $\Gamma_{2}$ of the domain and range respectively so that, up
to $\varepsilon$-homotopy, $L^{i}$ may be realized as a simplicial map sending
each $k$-simplex of $\Gamma_{1}$ onto a $k$-simplex of $\Gamma_{2}$. (See
Chapter 5 of \cite{RS}.) Adjust $\tau$ (rel $p_{i+1}$) so that it is an
embedded circle in general position with respect to $\Gamma_{2}$. Then
$\left(  L^{i}\right)  ^{-1}\left(  \tau\right)  $ is a closed $1$-manifold in
$\partial N_{i+2}\times(0,1)$. Let $\sigma$ be the component of $\left(
L^{i}\right)  ^{-1}\left(  \tau\right)  $ containing the point $\left(
p_{i+2},\frac{1}{4}\right)  $. Since $L^{i}$ takes a neighborhood of $\left(
p_{i+2},\frac{1}{4}\right)  $ homeomorphically onto a neighborhood of
$p_{i+1}$, and since no other points of $\sigma$ are taken near $p_{i+1}$ (use
condition v) from above), then $L$ restricts to a degree $\pm1$ map of
$\sigma$ onto $\tau$. Now the natural deformation retraction of $\partial
N_{i+2}\times\left[  0,1\right]  $ onto $\partial N_{i+2}\times\left\{
0\right\}  $ pushes $\sigma$ into $\partial N_{i+2}\times\left\{  0\right\}  $
while sliding $\left(  p_{i+2},\frac{1}{4}\right)  $ along the arc $\left\{
p_{i+2}\right\}  \times\left[  0,\frac{1}{4}\right]  $. Composing this push
with $L^{i}$ provides a homotopy of $\tau$ (within $N_{i}$) into $\partial
N_{i+2}$ whereby $p_{i+1}$ is slid along $\alpha_{i+1}$ to $p_{i+2}$.
\end{proof}

\begin{remark}
The reader may have noticed that a general principle at work in the proof of
Proposition \ref{2} is that ``degree $1$ maps between manifolds induce
surjections on fundamental groups''. Instead of applying this directly, we
used a constructive approach to finding the preimage of a loop. This allowed
us to handle orientable and non-orientable cases simultaneously. Proposition
\ref{3} is based on a similar general principle regarding homology groups and
degree $1$ maps. However, instead of a unified approach, we first obtain the
result for orientable manifolds by applying the general principle directly;
then we use the orientable result to extend to the non-orientable case. Those
who prefer this approach may use the proof of claim 1 from Proposition \ref{3}
as an outline to obtain an alternative proof of Proposition \ref{2} in the
case that $M^{n}$ is orientable.
\end{remark}

\begin{proposition}
\label{3}Let $M^{n}$ be a one ended, inward tame $n$-manifold with compact
boundary and let $R$ be a commutative ring with unity. Then $H_{j}\left(
\varepsilon\left(  M^{n}\right)  ;R\right)  $ is stable for all $i$.
\end{proposition}

For the sake of simplicity, we will first prove Proposition \ref{3} for
$R=\mathbb{Z}$. The more general result will then obtained by an application
of the universal coefficient theorem. Alternatively, one could do all of what
follows over an arbitrary coefficient ring. Before beginning the proof we
review some of the tools needed

Let $W$ be a compact connected orientable $n$-manifold with boundary. Assume
that $\partial W=P\cup Q$, where $P$ and $Q$ are disjoint, closed, $\left(
n-1\right)  $-dimensional submanifolds of $\partial W$. We do not require that
$P$ or $Q$ be connected or non-empty. Then Poincar\'{e} duality tells us that
the cap product with an orientation class $\left[  W\right]  $ induces
isomorphisms%
\[
H^{k}\left(  W,P;\mathbb{Z}\right)  \overset{\cap\;[W]}{\longrightarrow
}H_{n-k}\left(  W,Q;\mathbb{Z}\right)  \text{.}%
\]
If $W^{\prime}$ is another orientable $n$-manifold with $\partial W^{\prime
}=P^{\prime}\cup Q^{\prime}$, and $f:\left(  W,\partial W\right)
\rightarrow\left(  W^{\prime},\partial W^{\prime}\right)  $ is a map with
$f\left(  P\right)  \subset P^{\prime}$ and $f\left(  Q\right)  \subset
Q^{\prime}$, then the naturality of the cap product gives a commuting
diagram:
\begin{equation}%
\begin{array}
[c]{ccc}%
H^{k}\left(  W,P;\mathbb{Z}\right)  & \overset{\cap\;[W]}{\longrightarrow} &
H_{n-k}\left(  W,Q;\mathbb{Z}\right)  \text{.}\\
f^{\ast}\uparrow &  & \downarrow\,f_{\ast}\\
H^{k}\left(  W^{\prime},P^{\prime};\mathbb{Z}\right)  & \overset{\cap
\;f_{\ast}[W]}{\longrightarrow} & H_{n-k}\left(  W^{\prime},Q^{\prime
};\mathbb{Z}\right)
\end{array}
\tag{$\ddagger$}%
\end{equation}
If $f$ is of degree $\pm1$, then both horizontal homomorphisms are
isomorphisms, and hence $f_{\ast}$ is surjective.

For non-orientable manifolds, one may obtain duality isomorphisms and a
diagram like ($\ddagger$) by using $\mathbb{Z}_{2}$-coefficients. A more
powerful duality theorem and corresponding version of ($\ddagger$) for
non-orientable manifolds may be obtained by using ``twisted integer''
coefficients. This will be discussed after we handle the orientable case.

\begin{proof}
[Proof of Proposition \ref{3} (orientable case with
$\mathbb{Z}$-coefficients)]Let $M^{n}$ be orientable and let $\left\{
N_{i}\right\} _{i=0}^{\infty}$ be a sequence of neighborhoods of
infinity along with the embeddings, rays, base points, subspaces and
homotopies $\left\{ K^{i}\right\} _{i=0}^{\infty} $ described
earlier. For each $j\geq0$, $H_{j}\left( \varepsilon\left(
M^{n}\right) ;\mathbb{Z}\right) $ is represented by
\[
H_{j}\left(  N_{0};\mathbb{Z}\right)  \overset{\lambda_{1}}{\longleftarrow
}H_{j}\left(  N_{1};\mathbb{Z}\right)  \overset{\lambda_{2}}{\longleftarrow
}H_{j}\left(  N_{2};\mathbb{Z}\right)  \overset{\lambda_{3}}{\longleftarrow
}\cdots
\]
where all bonding maps are induced by inclusion.

Since each $N_{i}$ is connected, $H_{0}\left(  \varepsilon\left(
M^{n}\right)  ;\mathbb{Z}\right)  $ is pro-equivalent to $$\mathbb{Z}%
\overset{\cong}{\leftarrow}\mathbb{Z\overset{\cong}{\leftarrow}Z\overset
{\cong}{\leftarrow}\cdots}$$
and thus, is stable. Let $j\geq1$ be fixed.\medskip

\noindent\textbf{Claim 1}\qua\textbf{\ }$H_{j}\left(  \varepsilon\left(
M^{n}\right)  ;\mathbb{Z}\right)  $ {\sl is semistable.}\medskip

We will show that, for each $\left[  \alpha\right]  \in H_{j}\left(
N_{i+1}\right)  $, there is a $\left[  \alpha^{\prime}\right]  \in
H_{j}\left(  N_{i+2}\right)  $ such that $\alpha$ is homologous to
$\alpha^{\prime}$ in $N_{i}$. Thus, $im\left(  \lambda_{i+1}\right)
\overset{\lambda_{i+1}}{\longleftarrow}im\left(  \lambda_{i+2}\right)  $ is surjective.

We may assume that $\alpha$ is supported in $A_{i+1}$. We abuse notation
slightly and write $\left[  \alpha\right]  \in H_{j}\left(  A_{i+1}%
;\mathbb{Z}\right)  $. Let $L^{i}:\partial N_{i+2}\times\left[  0,1\right]
\rightarrow N_{i}$ be the restriction of $K^{i}$. Note that $L^{i}\left(
\partial N_{i+2}\times\left\{  1\right\}  \right)  \subset A_{i}-\partial
N_{i+1}$. By PL transversality theory (see \cite{RS1} or Section II.4 of
\cite{BRS}), we may---after a small adjustment that does not alter $L^{i}$ on
$(\partial N_{i+2}\times\left\{  0,1\right\}  )\cup(D_{i}\times\left[
0,1\right]  )$---assume that that $C_{i+1}\equiv(L^{i})^{-1}\left(
A_{i+1}\right)  $ is an $n$-manifold with boundary\footnote{Instead of using
transversality theory, we could simply use the radial structure of regular
neighborhoods to alter $L^{i}$ in \ a thin regular neighborhood of
$(L^{i})^{-1}\left(  A_{i}\cup N_{i+2}\right)  $. Using this approach, we
``fatten'' the preimage of $A_{i}\cup N_{i+2}$ to a codimension $0$
submanifold, thus ensuring that $(L^{i})^{-1}\left(  A_{i+1}\right)  $ is an
$n$-manifold with boundary.}. Let $C_{i+1}^{\ast}$ be the component of
$C_{i+1}$ that contains $D_{i}^{\prime}\times\left[  0,\frac{1}{4}\right]  $.
Then $L^{i}$ takes $\partial C_{i+1}^{\ast}$ into $\partial A_{i+1}$ and,
provided our adjustment to $L^{i}$ was sufficiently small, $L^{i}$ is still a
homeomorphism \emph{over} $T_{0}^{\prime}\cap A_{i+1}$. By the local
characterization of degree, $L^{i}\mid_{C_{i+1}^{\ast}}:\left(  C_{i+1}^{\ast
},\partial C_{i+1}^{\ast}\right)  \rightarrow\left(  A_{i+1},\partial
A_{i+1}\right)  $ is a degree $1$ map. Thus, by an application of ($\ddagger
$), $\left[  \alpha\right]  $ has a preimage $\left[  \beta\right]  \in
H_{j}\left(  C_{i+1}^{\ast};\mathbb{Z}\right)  $. Now $C_{i+1}^{\ast}%
\subset\partial N_{i+2}\times\left[  0,1\right]  $, and within the larger
space, $\beta$ is homologous to a cycle $\beta^{\prime}$ supported in
$\partial N_{i+2}\times\left\{  0\right\}  $. Since $L^{i}$ takes $\partial
N_{i+2}\times\left[  0,1\right]  $ into $N_{i}$, it follows that $\alpha$ is
homologous to $\alpha^{\prime}\equiv L^{i}\left(  \beta^{\prime}\right)
\subset\partial N_{i+2}$ in $N_{i}$.\medskip

\noindent\textbf{Claim 2}\qua $H_{j}\left(  \varepsilon\left(  M^{n}\right)
;\mathbb{Z}\right)  $ {\sl is pro-monomorphic.\medskip}

We'll show that $im\left(  \lambda_{i+2}\right)  \overset{\lambda_{i+1}%
}{\longleftarrow}im\left(  \lambda_{i+3}\right)  $ is injective, for all
$i\geq0$. It suffices to show that each $j$-cycle $\alpha$ in $N_{i+3}$ that
bounds a $\left(  j+1\right)  $-chain $\gamma$ in $N_{i+1}$, bounds a $\left(
j+1\right)  $-chain in $N_{i+2}$. Let $[\gamma^{\prime}]$ be a preimage of
$\left[  \gamma\right]  $ under the excision isomorphism%
\[
H_{j+1}\left(  A_{i+1}\cup A_{i+2},\partial N_{i+3};\mathbb{Z}\right)
\rightarrow H_{j+1}\left(  N_{i+1},N_{i+3};\mathbb{Z}\right)  \text{.}%
\]
Then $\alpha^{\prime}\equiv\partial\gamma^{\prime}$ is homologous to $\alpha$
in $N_{i+3}$, so it suffices to show that $\alpha^{\prime}$ bounds in
$N_{i+2.}$

By passing to a subsequence if necessary, we may assume that the image of
$\partial N_{i+2}\times\left[  0,1\right]  $ under $K^{i}$ lies in $A_{i}\cup
A_{i+1}\cup A_{i+2}-U$, where $U$ is a collar neighborhood of $\partial
N_{i+3}$ in $A_{i+2}$. Then define
\[
f:(\partial N_{i+2}\times\left[  0,1\right]  )\cup A_{i+2}\rightarrow
A_{i}\cup A_{i+1}\cup A_{i+2}%
\]
to be $K^{i}$ on $\partial N_{i+2}\times\left[  0,1\right]  $ and the identity
on $A_{i+2}$. Arguing as in the proof of Claim 1, we may---without changing
the map on $A_{i+2}$\ ---make a small adjustment to $f$ so that $C\equiv
f^{-1}\left(  A_{i+1}\cup A_{i+2}\right)  $ is an $n$-manifold with boundary.
Let $C^{\ast}$ be the component that contains $A_{i+2}$. Then $f$ takes
$\partial N_{i+3}$ onto $\partial N_{i+3}$, and $P\equiv\partial C^{\ast
}-\partial N_{i+3}$ to $\partial N_{i+1}$. Provided our adjustment was
sufficiently small, $f$ is a homeomorphism \emph{over} $U$, so $f:\left(
C^{\ast},\partial C^{\ast}\right)  \rightarrow\left(  A_{i+1}\cup
A_{i+2},\partial N_{i+1}\cup\partial N_{i+3}\right)  $ is a degree $1$ map
Applying ($\ddagger$) to this situation we obtain a surjection
\[
H_{j+1}\left(  C,\partial N_{i+3};\mathbb{Z}\right)  \twoheadrightarrow
H_{j+1}\left(  A_{i+1}\cup A_{i+2},\partial N_{i+3};\mathbb{Z}\right)
\text{.}%
\]
Let $[\eta]$ be a preimage of $[\alpha^{\prime}]$. Utilizing the product
structure on $\partial N_{i+2}\times\left[  0,1\right]  $, we may retract
$C^{\ast}$ onto $A_{i+2}$. The image $\eta^{\prime}$ of $\eta$ under this
retraction is a relative $\left(  j+1\right)  $-cycle in $\left(
A_{i+2},\partial N_{i+3}\right)  $ with $\partial\eta^{\prime}=\partial\eta$.
Thus, $\partial\eta^{\prime}$ is homologous to $\partial\gamma^{\prime}%
=\alpha^{\prime}$, so $\alpha^{\prime}$ bounds in $A_{i+2}\subset N_{i+2}$ as desired.
\end{proof}

Before proceeding with the proof of the non-orientable case, we discuss some
necessary background. The proof just presented already works for
non-orientable manifolds if we replace the coefficient ring $\mathbb{Z}$ with
$\mathbb{Z}_{2}$. To obtain the result for $\mathbb{Z}$-coefficients (and
ultimately an arbitrary coefficient ring), we will utilize homology with
twisted integer coefficients, which we will denote by $\widetilde{\mathbb{Z}}%
$. The key here is that, even for a \emph{non-orientable} compact $n$-manifold
with boundary, $H_{n}\left(  W,\partial W;\widetilde{\mathbb{Z}}\right)
\cong\mathbb{Z}$. Thus, we have an orientation class $\left[  W\right]  $ and
it may be used to obtain a duality isomorphism---where homology is now taken
with twisted integer coefficients. Furthermore, if a map $f:\left(  W,\partial
W\right)  \rightarrow\left(  W^{\prime},\partial W^{\prime}\right)  $ is
\emph{orientation true }(meaning that $f$ takes orientation reversing loops to
orientation reversing loops and orientation preserving loops to orientation
preserving loops), then we have a well defined notion of $\deg\left(
f\right)  \in\mathbb{Z}$. These versions of duality and degree yield an
analogous version of diagram ($\ddagger$), which tells us that degree $\pm1$
maps (appropriately defined) between compact (possibly non-orientable)
manifolds with boundary induce surjections on homology with $\widetilde
{\mathbb{Z}}$-coefficients. See section 3.H of \cite{Ha} or Chapter 2 of
\cite{Wl2} discussions of homology with coefficients in $\widetilde
{\mathbb{Z}}$, and \cite{Ol} for a discussion of degree of a map between
non-orientable manifolds.

As with the traditional definition of degree, this generalized version can be
detected locally. In particular, an orientation true map $f:(W,\partial
W)\rightarrow(W^{\prime},\partial W^{\prime})$ that is a homeomorphism over
some open subset of $W^{\prime}$ has degree $\pm1$. See \cite[3.8]{Ol}.

For non-orientable $W$, let $p:\widehat{W}\rightarrow W$ be the orientable
double covering projection. Then there is a long exact sequence:%
\[
\cdots\rightarrow H_{k}\left(  W;\widetilde{\mathbb{Z}}\right)  \rightarrow
H_{k}\left(  \widehat{W};\mathbb{Z}\right)  \overset{p_{\ast}}{\rightarrow
}H_{k}\left(  W;\mathbb{Z}\right)  \rightarrow H_{k-1}\left(  W;\widetilde
{\mathbb{Z}}\right)  \rightarrow\cdots
\]
This sequence is natural with respect to orientation true mappings
$f:W\rightarrow W^{\prime}$. See section 3.H of \cite{Ha} for a discussion of
this sequence.

\begin{proof}
[Proof of Proposition \ref{3} (non-orientable case with 
$\mathbb{Z}$-coefficients)]\noindent\hbox{}\hfill\break
Let $M^{n}$ be one ended, inward tame, and have compact boundary. If $M^{n}$
contains an orientable neighborhood of infinity, we can simply disregard its
complement and apply the orientable case. Hence, we assume that $\left\{
N_{i}\right\}  _{i=0}^{\infty}$ is a nested cofinal sequence of $0$%
-neighborhoods of infinity, each of which is non-orientable.

The first step in this proof is to observe that, if we use homology with
$\widetilde{\mathbb{Z}}$-coefficients, the proof used in the orientable case
is still valid. A few points are worth noting. First, the inclusion maps
$N_{i}\hookrightarrow N_{i+1}$ are clearly orientation true. Similarly, since
each $\partial N_{i}$ is bicollared in $M^{n}$, orientation reversing
[preserving] loops in $\partial N_{i}$ are orientation reversing [preserving]
in $M^{n}$. Hence, the maps $L^{i}:\partial N_{i+2}\times\left[  0,1\right]
\rightarrow N_{i}$ (and restrictions to codimension $0$ submanifolds) are also
orientation true. With this, and the additional ingredients discussed above,
we see that $H_{j}\left(  \varepsilon(M^{n});\widetilde{\mathbb{Z}}\right)  $
is stable for all $j$.

The second step is to consider the orientable double covering projection
$p:\widehat{M}^{n}\rightarrow M^{n}$. For each $i$, $\widehat{N_{i}}%
=p^{-1}(N_{i})$ is the orientable double cover of $N_{i}$, and thus, is
connected. It follows that $\widehat{M}^{n}$ is one ended, with $\left\{
\widehat{N_{i}}\right\}  _{i=0}^{\infty}$ a cofinal sequence of $0$%
-neighborhoods of infinity. Furthermore, taming homotopies for $M^{n}$ may be
lifted to obtain taming homotopies for $\widehat{M}^{n}$, so $\widehat{M}^{n}$
is inward tame. It follows from the orientable case that $H_{j}\left(
\varepsilon(\widehat{M}^{n});\mathbb{Z}\right)  $ is stable for all $j$.

Next we apply the long exact discussed above to each covering projection
$p_{i}:\widehat{N_{i}}\rightarrow N$. Together with naturality, this yields a
long exact sequence of inverse sequences:%
\vspace{-5mm}
{\small
\[%
\begin{array}
[c]{ccccccccc}
 & \vdots &  & \vdots &  & \vdots &  & \vdots &   \\
 & \downarrow &  & \downarrow &  & \downarrow &  & \downarrow &   \\
\cdots \rightarrow & H_{k}\left(  N_{3};\widetilde{\mathbb{Z}}\right)  &
\rightarrow & H_{k}\left(  \widehat{N}_{3};\mathbb{Z}\right)  & \rightarrow &
H_{k}\left(  N_{3};\mathbb{Z}\right)  & \rightarrow & H_{k-1}\left(
N_{3};\widetilde{\mathbb{Z}}\right)  & \rightarrow \cdots\\
 & \downarrow &  & \downarrow &  & \downarrow &  & \downarrow &   \\
\cdots \rightarrow & H_{k}\left(  N_{2};\widetilde{\mathbb{Z}}\right)  &
\rightarrow & H_{k}\left(  \widehat{N}_{2};\mathbb{Z}\right)  & \rightarrow &
H_{k}\left(  N_{2};\mathbb{Z}\right)  & \rightarrow & H_{k-1}\left(
N_{2};\widetilde{\mathbb{Z}}\right)  & \rightarrow \cdots\\
 & \downarrow &  & \downarrow &  & \downarrow &  & \downarrow  & \\
\cdots \rightarrow & H_{k}\left(  N_{1};\widetilde{\mathbb{Z}}\right)  &
\rightarrow & H_{k}\left(  \widehat{N}_{1};\mathbb{Z}\right)  & \rightarrow &
H_{k}\left(  N_{1};\mathbb{Z}\right)  & \rightarrow & H_{k-1}\left(
N_{1};\widetilde{\mathbb{Z}}\right)  & \rightarrow \cdots
\end{array}
\]}%
We may now apply Lemma \ref{dydak} to conclude that $H_{j}\left(
\varepsilon(M^{n});\mathbb{Z}\right)  $ is stable for all $j$.
\end{proof}

Lastly, we generalize the above to the case of an arbitrary coefficient ring.

\begin{proof}
[Proof of Proposition \ref{3} ($R$-coefficients)]Now let $R$ be ring with unity. By
applying the Universal Coefficient Theorem for homology (see \cite[Cor.
3.A.4]{Ha}) to obtain each row, we may get (for each $j$) the following
diagram:%
\[%
\begin{array}
[c]{ccccccccc}%
\vdots &  & \vdots &  & \vdots &  & \vdots &  & \vdots\\
\downarrow &  & \downarrow &  & \downarrow &  & \downarrow &  & \downarrow\\
0 & \rightarrow & H_{j}\left(  N_{3};\mathbb{Z}\right)  \otimes R &
\rightarrow & H_{j}\left(  N_{3};R\right)  & \rightarrow & Tor(H_{j-1}\left(
N_{3};\mathbb{Z}\right)  ,R) & \rightarrow & 0\\
\downarrow &  & \downarrow &  & \downarrow &  & \downarrow &  & \downarrow\\
0 & \rightarrow & H_{j}\left(  N_{2};\mathbb{Z}\right)  \otimes R &
\rightarrow & H_{j}\left(  N_{2};R\right)  & \rightarrow & Tor(H_{j-1}\left(
N_{2};\mathbb{Z}\right)  ,R) & \rightarrow & 0\\
\downarrow &  & \downarrow &  & \downarrow &  & \downarrow &  & \downarrow\\
0 & \rightarrow & H_{j}\left(  N_{1};\mathbb{Z}\right)  \otimes R &
\rightarrow & H_{j}\left(  N_{1};R\right)  & \rightarrow & Tor(H_{j-1}\left(
N_{1};\mathbb{Z}\right)  ,R) & \rightarrow & 0
\end{array}
\]

The second and fourth columns are stable by the $\mathbb{Z}$-coefficient case,
so an application of Lemma \ref{dydak} yields stability of $H_{j}\left(
\varepsilon(M^{n});R\right)  .$
\end{proof}

\begin{remark}
A variation on the above can be used to show that, for one ended manifolds
with compact boundary, inward tameness plus $\pi_{1}$-stability implies
$\pi_{2}$-stability. To do this, begin with a cofinal sequence $\left\{
N_{i}\right\}  $ of (strong) $1$-neighborhoods of infinity---see Theorem 4 of
\cite{Gu1}. Then show that the inverse sequence%
\[
H_{2}\left(  \widetilde{N}_{0};\mathbb{Z}\right)  \leftarrow H_{2}\left(
\widetilde{N}_{1};\mathbb{Z}\right)  \leftarrow H_{2}\left(  \widetilde{N}%
_{2};\mathbb{Z}\right)  \leftarrow\cdots
\]
is stable, where each $\widetilde{N}_{i}$ is the universal cover of $N_{i}$.
This will require Poincar\'{e} duality for noncompact manifolds; otherwise,
the proof simply mimics the proof of Prop. \ref{3}. It follows from the
Hurewicz theorem that
\[
\pi_{2}\left(  \widetilde{N}_{0},\tilde{p}_{0}\right)  \leftarrow\pi
_{2}\left(  \widetilde{N}_{1},\tilde{p}_{1}\right)  \leftarrow\pi_{2}\left(
\widetilde{N}_{2},\tilde{p}_{2}\right)  \leftarrow\cdots
\]
is stable, and hence, so is $\pi_{2}\left(  \varepsilon\left(  M^{n}\right)
\right)  $. As an application of this observation, one may deduce the main
result of Siebenmann's thesis as a direct corollary of Theorem \ref{pseudo}%
---provided $n\geq7$.
\end{remark}

\section{Proof of Theorem \ref{perfect}\label{construction}}

In this section we will construct (for each $n\geq6$) a one ended open
$n$-manifold $M_{\ast}^{n}$ in which all clean neighborhoods of infinity have
finite homotopy type, yet $\pi_{1}(\varepsilon(M_{\ast}^{n}))$ is not
perfectly semistable. Hence $M_{\ast}^{n}$ satisfies conditions (1) and (3) of
Theorem \ref{pseudo}, but is not pseudo-collarable.

In the first portion of this section we present the necessary group theory on
which the examples rely. In the next portion, we give a detailed construction
of the examples and verify the desired properties.

\subsection{Group Theory}

We assume the reader is familiar with the basic notions of group presentations
in terms of generators and relators. We use the HNN-extension as our basic
building block. A more thorough discussion of HNN-extensions may be found in
\cite{KS} or \cite{Co}.

Before beginning, we describe the algebraic goal of this section. We wish to
construct a special inverse sequence of finitely presented groups that is
semistable, but not perfectly semistable. Later this sequence will be realized
as the fundamental group at infinity of a carefully constructed open manifold.
The following lemma indicates the strategy that will be used.

\begin{lemma}
\label{L3}Let
\[
G_{0}\overset{\psi_{1}}{\longleftarrow}G_{1}\overset{\psi_{2}}{\longleftarrow
}G_{2}\overset{\psi_{3}}{\longleftarrow}G_{3}\overset{\psi_{4}}{\longleftarrow
}\cdots
\]
be an inverse sequence of groups with surjective but non-injective bonding
homomorphisms. Suppose further that no $G_{i}$ contains a non-trivial perfect
subgroup. Then this inverse sequence is not perfectly semistable.

\begin{proof}
It is easy to see that this system is semistable but not stable. Assume that
it is perfectly semistable. Then---after passing to a subsequence,
relabelling, and applying Lemma \ref{L1}---we may assume the existence of a
diagram:
\[%
\begin{array}
[c]{ccccccc}%
G_{0} & \overset{\psi_{1}}{\twoheadleftarrow} & G_{1} & \overset{\psi_{2}%
}{\twoheadleftarrow} & G_{2} & \overset{\psi_{3}}{\twoheadleftarrow} &
\cdots\\
& \underset{f_{0}\quad}{\nwarrow}\quad\underset{\quad g_{0}}{\swarrow} &  &
\underset{f_{1}\quad}{\nwarrow}\quad\underset{\quad g_{1}}{\swarrow} &  &
\underset{f_{2}\quad}{\nwarrow}\quad\underset{\quad g_{2}}{\swarrow} & \\
& H_{0} & \overset{\mu_{1}}{\twoheadleftarrow} & H_{1} & \overset{\mu_{2}%
}{\twoheadleftarrow} & H_{2} & \cdots
\end{array}
\]
where each $\mu_{i}$ has a perfect kernel.

By the commutativity of the diagram, all of the $f_{i}$'s and $g_{i}$'s are
surjections. Moreover, Lemma \ref{L2} implies that $f_{i}\left(  \ker\mu
_{i}\right)  =\left\{  1\right\}  $, for all $i\geq1$. The combination of
these facts tells us that each $g_{i}$ is an isomorphism. Since the $G_{i}$'s
contain no nontrivial perfect subgroups, then neither do the $H_{i}$'s. But
then each $\mu_{i}$ is an isomorphism, contradicting the non-stability of our
original sequence.
\end{proof}
\end{lemma}

\begin{remark}
Satisfying the hypotheses of Lemma \ref{L3} by itself is not difficult. For
example, since abelian groups contain no nontrivial perfect subgroups,
examples such as
\[
\mathbb{Z\twoheadleftarrow Z\oplus Z\twoheadleftarrow Z\oplus Z\oplus
Z\twoheadleftarrow\cdots}%
\]
apply. However, Theorem \ref{semistable} tells us that this inverse sequence
cannot occur as the fundamental group at infinity of an inward tame open
manifold. Indeed, any appropriate inverse sequence should, at least, have the
property that abelianizing each term yields a stable sequence. Thus, our task
of constructing an appropriate ``realizable'' inverse sequence is rather delicate.
\end{remark}

Let $K$ be a group with presentation $\left\langle gen(K)|rel(K)\right\rangle
$ and $\{\phi_{i}\}$ a collection of monomorphisms $\phi_{i}:L_{i}\rightarrow
K$ from subgroups $\{L_{i}\}$ of $K $ into $K$. We define the group%
\[
G=\left\langle gen(K),t_{1},t_{2},\cdots\mid rel(K),R_{1},R_{2},\cdots
\right\rangle
\]
where each $R_{i}$ is the collection of relations $\left\{  t_{i}l_{ij}%
t_{i}^{-1}=\phi_{i}\left(  l_{ij}\right)  \text{ for all }l_{ij}\in
L_{i}\right\}  $. We call $G$ the \emph{HNN } \emph{group} with \emph{base K,
associated subgroups} $\left\{  L_{i},\phi_{i}\left(  L_{i}\right)  \right\}
$\emph{, \/} and \emph{free part\/} the group generated by $\left\{
t_{1},t_{2},\cdots\right\}  $. We assume the basic properties of HNN
groups---such as the fact that the base group naturally embeds in the HNN
group. This and other basic structure theorems for subgroups of HNN-extensions
have existed for a long time and appear within many sources. Most important
for our purposes is the following which we have tailored to meet our specific needs.

\begin{theorem}
[{see \cite[Theorem 6]{KS}}]\label{subgroups}Let $G$ be the HNN group above. If
$H$ is a subgroup of $G$ having trivial intersection with the conjugates of
each $L_{i}$, then $H$ is the free product of a free group with the
intersections of $H$ with certain conjugates of $K$.
\end{theorem}

Let $a$ and $b$ be group elements. We denote by $[a,b]$ the
\emph{commutator\/} of $a$ and $b$, ie, $\left[  a,b\right]  =a^{-1}%
b^{-1}ab$. Let $S$ be a subset of elements of a group $G$. We denote by
$\left\{  s_{1},s_{2},\cdots;G\right\}  $ the \emph{subgroup of }%
$G$\emph{\ generated by }$S$ where $S=\left\{  s_{1},s_{2},\cdots\right\}  $.
If $S$ and $G$ are as above, then we denote by $ncl\left\{  s_{1},s_{2}%
,\cdots;G\right\}  $ the \emph{normal closure of }$S$\emph{\ in }$G$, ie,
the smallest normal subgroup of $G$ containing $S.$

We are now ready to construct the desired inverse sequence. Let $G_{0}%
=\left\langle a_{0}\right\rangle $ be the free group on one generator. Of
course, $G_{0}$ is just $\mathbb{Z}$ written multiplicatively. For $j\geq1$,
let
\[
G_{j}=\left\langle a_{0},a_{1},\cdots,a_{j}\mid a_{1}=\left[  a_{1}%
,a_{0}\right]  ,a_{2}=\left[  a_{2},a_{1}\right]  ,\cdots,a_{j}=\left[
a_{j},a_{j-1}\right]  \right\rangle
\]
This presentation emphasizes that each $a_{i}$ ($i\geq1$) is a commutator.
(Hence, each $G_{j}$ abelianizes to $\mathbb{Z}$.) We abuse notation slightly
and do not distinguish between the element $a_{i}\in G_{j-1}$ and $a_{i}\in
G_{j}.$ Let $j\geq1$; another useful presentation of $G_{j}$ is
\[
G_{j}=\left\langle a_{0},a_{1},\cdots,a_{j}\mid a_{0}a_{1}^{2}a_{0}^{-1}%
=a_{1},a_{1}a_{2}^{2}a_{1}^{-1}=a_{2},\cdots,a_{j-1}a_{j}^{2}a_{j-1}%
^{-1}=a_{j}\right\rangle
\]
Now, $G_{j}$ can be put in the form of an HNN group. In particular,
\[
G_{j}=\left\langle gen(K),t_{1}\mid rel(K),R_{1}\right\rangle
\]
where
\[
K=\left\langle a_{1},a_{2},\cdots,a_{j}\mid a_{1}a_{2}^{2}a_{1}^{-1}%
=a_{2},a_{2}a_{3}^{2}a_{2}^{-1}=a_{3},\cdots,a_{j-1}a_{j}^{2}a_{j-1}%
^{-1}=a_{j}\right\rangle ,
\]
$t_{1}=a_{0}$, $L_{1}=\left\{  a_{1}^{2};G_{j}\right\}  $, $\phi_{1}\left(
a_{1}^{2}\right)  =a_{1}$, and $R_{1}$ is given by $a_{0}a_{1}^{2}a_{0}%
^{-1}=a_{1}$. The base group, $K$, is obviously isomorphic to $G_{j-1}$ with
that isomorphism taking $a_{i}$ to $a_{i-1}$.

Define $\psi_{j}:G_{j}\rightarrow G_{j-1}$ by sending $a_{i}$ to $a_{i}$ for
$1\leq i\leq j-1$, and $a_{j}$ to $1$. By inspection $\psi_{j}$ is a
surjective homomorphism. Our goal is to prove:

\begin{theorem}
\label{no-perfects}In the setting described above, the group $G_{j}$ has no
non-trivial perfect subgroups.

\begin{proof}
Our proof is by induction.\medskip

\noindent\textbf{Case $j=0$}\quad$G_{0}=\left\langle a_{0}%
\right\rangle $ is an abelian group so that all commutators in $G_{0}$ are
trivial. Thus, $[H,H]=1$ for any subgroup $H$ of $G_{0}$. Hence, $H=1$ is the
only perfect subgroup of $G_{0}$.\medskip

\noindent\textbf{Case $j=1$}\quad Consider $G_{1}$ and $\psi
_{1}:G_{1}\rightarrow G_{0}$. $\psi_{1}:\left\langle a_{0},a_{1}|a_{0}%
a_{1}^{2}a_{0}^{-1}=a_{1}\right\rangle $$\rightarrow\left\langle
a_{0}\right\rangle $. We pause to observe for later use that $G_{1}$ is an HNN
group with base group $K=\left\{  a_{1};G_{1}\right\}  $. Since $K$ embeds in
$G_{1}$, then $a_{1}$ has infinite order in $G_{1}.$ Now, $G_{1}$ is one of
the well-known Baumslag-Solitar groups. Its commutator subgroup, $[G_{1}%
,G_{1}]$, is precisely equal to $ker\left(  \psi_{1}\right)  $. The
substitution $b_{k}=a_{0}^{-k}a_{1}a_{0}^{k}$ along with the relations
\[
b_{k}=a_{0}^{-k}a_{1}a_{0}^{k}=a_{0}^{-(k-1)}\big(     a_{0}^{-1}a_{1}%
a_{0}^{1}\big)     a_{0}^{k-1}=a^{-(k-1)}a_{1}^{2}a^{k-1}=\big(
a^{-(k-1)}a_{1}a^{k-1}\big)     ^{2}%
\]
give $ker\left(  \psi_{1}\right)  $ a presentation:
\[
\left\langle b_{k}\mid b_{k}=b_{k-1}^{2},-\infty<j<\infty\right\rangle
\]
So, $ker\left(  \psi_{1}\right)  $ is locally cyclic (every finitely generated
subgroup is contained in a cyclic subgroup). In particular, it is abelian and
contains no non-trivial perfect subgroups. Now, suppose $P$ is a perfect
subgroup of $G_{1}$. Then, by Lemma \ref{L2}, $\psi_{1}\left(  P\right)  $ is
a perfect subgroup of $G_{0}$. By the case ($j=0$), $\psi_{1}\left(  P\right)
=\left\{  1\right\}  $, so $P\subset ker\left(  \psi_{1}\right)  $. But, we
just observed, then, that $P$ must be trivial.\medskip

\noindent\textbf{Inductive Step}\qua We assume that $G_{j}$ contains no
non-trivial perfect subgroups for $1\leq j\leq k-1$ and prove that $G_{k}$ has
this same property. To this end, let $P$ be a perfect subgroup of $G_{k}$.
Then, $\psi_{k}\left(  P\right)  $ is a perfect subgroup of $G_{k-1}$. By
induction, $\psi_{k}(P)=1$. Thus, $P\subset ker\left(  \psi_{k}\right)  $.

As shown above, $G_{k}$, is an HNN-extension with base group $K$ where
\begin{align*}
K  &  =\left\langle a_{1},a_{2},\cdots,a_{k}\mid a_{1}a_{2}^{2}a_{1}%
^{-1}=a_{2},a_{2}a_{3}^{2}a_{2}^{-1}=a_{3},\cdots,a_{k-1}a_{k}^{2}a_{k-1}%
^{-1}=a_{k}\right\rangle \\
&  \cong G_{k-1}%
\end{align*}
By the inductive hypothesis, $K$ has no perfect subgroups. Moreover, $a_{1}\in
K$ still has infinite order in both $K$ (by induction) and $G_{k}$ (since $K$
embeds in $G_{k}$). Moreover, the HNN group, $G_{k}$, has the single
associated cyclic subgroup, $L=\{a_{1}^{2};G_{k}\}$, with conjugation relation
$a_{0}a_{1}^{2}a_{0}^{-1}=a_{1}$. By the definition of $\psi_{k}%
:G_{k}\rightarrow G_{k-1}$ it is clear that $ker\left(  \psi_{k}\right)
=ncl\big\{     a_{k};G_{k}\big\}     $.\medskip

\noindent\textbf{Claim}\qua {\sl No conjugate of }$L${\ \sl non-trivially
intersects }$ncl\big\{     a_{k};G_{k}\big\}     \medskip$

\noindent\textbf{Proof of Claim}\qua 
If the claim is false, then $L$ itself must
non-trivially intersect the normal subgroup, $ncl\left\{  a_{k};G_{k}\right\}
$. This means that $a_{1}^{2m}\in ncl\left\{  a_{k};G_{k}\right\}  =ker\left(
\psi_{k}\right)  $ for some integer $m>0$. Since $k\geq2$, then $\psi
_{k}\left(  a_{1}^{2m}\right)  =\psi_{k}\left(  a_{1}\right)  ^{2m}=a_{1}%
^{2m}=1 $ in $G_{k-1}$, ie, $a_{1}$ has finite order in $G_{k-1} $. This
contradicts our observations above, thus proving the claim.\medskip

We continue with the proof of Theorem \ref{no-perfects}. Recall that $P$ is a
perfect subgroup of $ker\left(  \psi_{k}\right)  $. It must also enjoy the
property of trivial intersection with each conjugate of $L$. We now apply
Theorem \ref{subgroups} to the subgroup $P$ to conclude that $P$ is a free
product where each factor is either free or equal to $P\cap gKg^{-1}$ for some
$g\in G_{k}$.

Now, $P$ projects naturally onto each of these factors so each factor is
perfect. However, non-trivial free groups are not perfect. Moreover, by
induction, $K$ (or equivalently $gKg^{-1}$) contains no non-trivial perfect
subgroups. Thus, any subgroup, $P\cap gKg^{-1}$, is trivial. Consequently, $P
$ must be trivial.
\end{proof}
\end{theorem}

\subsection{Construction of $M_{\ast}^{n}$}

The goal of this section is to construct a one ended open $n$-manifold
$M_{\ast}^{n}$ ($n\geq6$) with fundamental group system at infinity equivalent
to the inverse sequence
\begin{equation}
G_{0}\overset{\psi_{1}}{\longleftarrow}G_{1}\overset{\psi_{2}}{\longleftarrow
}G_{2}\overset{\psi_{3}}{\longleftarrow}G_{3}\overset{\psi_{4}}{\longleftarrow
}\cdots\tag{$\dagger{}\dagger$}%
\end{equation}
produced above. More importantly, this will be done in such a way that clean
neighborhoods of infinity in $M_{\ast}^{n}$ have finite homotopy
type---thereby proving Theorem \ref{perfect}. Familiarity with the basics of
handle theory, as can be found in Chapter 6 of \cite{RS}, is assumed
throughout the construction.

The key to producing $M_{\ast}^{n}$ will be a careful construction of a
sequence $$\left\{  \left(  A_{i},\Gamma_{i},\Gamma_{i+1}\right)  \right\}
_{i=0}^{\infty}$$ of compact $n$-dimensional cobordisms satisfying the
following properties:\smallskip

\begin{enumerate}
\item[a)] The left-hand boundary $\Gamma_{0}$ of $A_{0}$ is $S^{n-2}\times
S^{1},$ and (as indicated by the notation), for all $i\geq1$ the left-hand
boundary of $A_{i}$ is equal to the right-hand boundary of $A_{i-1}$. In
particular, $A_{i-1}\cap A_{i}=\Gamma_{i}$.

\item[b)] For all $i\geq0$, $\pi_{1}(\Gamma_{i},p_{i})\cong G_{i}$ and
$\Gamma_{i}\hookrightarrow A_{i}$ induces a $\pi_{1}$-isomorphism.

\item[c)] The isomorphisms between $\pi_{1}(\Gamma_{i},p_{i})$ and $G_{i}$ may
be chosen so that we have a commutative diagram:%
\[%
\begin{array}
[c]{ccc}%
G_{i} & \overset{\psi_{i+1}}{\longleftarrow} & G_{i+1}\\
\downarrow\cong &  & \downarrow\cong\\
\pi_{1}(\Gamma_{i},p_{i}) & \overset{\mu_{i+1}}{\longleftarrow} & \pi
_{1}(\Gamma_{i+1},p_{i+1})
\end{array}
\]
Here $\mu_{i+1}$ is the composition of homomorphisms%
\[
\pi_{1}\left(  \Gamma_{i},p_{i}\right)  \overset{\cong}{\leftrightarrow}%
\pi_{1}\left(  A_{i},p_{i}\right)  \overset{\widehat{\alpha}_{i}}{\leftarrow
}\pi_{1}\left(  A_{i},p_{i+1}\right)  \overset{j_{\#}}{\longleftarrow}\pi
_{1}\left(  \Gamma_{i+1},p_{i+1}\right)
\]
where $j_{\#}$ is induced by inclusion, the middle map is a ``change of base
points isomorphism'' with respect to a path $\alpha_{i}$ in $A_{i}$ between
$p_{i}$ and $p_{i+1}$, and the left-most isomorphism is provided by property b).\smallskip
\end{enumerate}

\noindent We will let
\[
M_{\ast}^{n}=(S^{n-2}\times B^{2})\cup A_{0}\cup A_{1}\cup A_{2}\cup\cdots
\]
where $S^{n-2}\times B^{2}$ is glued to $A_{0}$ along $\Gamma_{0}%
=S^{n-2}\times S^{1}$. Then for each $i\geq0,$
\[
N_{i}=A_{i}\cup A_{i+1}\cup A_{i+2}\cup\cdots
\]
is a clean connected neighborhood of infinity. Moreover, by properties b) and
c) and repeated application of the Seifert-VanKampen theorem, the inverse
sequence
\[
\pi_{1}\left(  N_{0},p_{0}\right)  \overset{\lambda_{1}}{\longleftarrow}%
\pi_{1}\left(  N_{1},p_{1}\right)  \overset{\lambda_{2}}{\longleftarrow}%
\pi_{1}\left(  N_{2},p_{2}\right)  \overset{\lambda_{3}}{\longleftarrow}\cdots
\]
is isomorphic to ($\dagger{}\dagger$).

Finally, we will need to show that clean neighborhoods of infinity in
$M_{\ast}^{n}$ have finite homotopy type. This can be done only after the
specifics of the construction are revealed.\medskip

\noindent\textbf{Step 0}\qua {\sl Construction of }$\left(  A_{0},\Gamma
_{0},\Gamma_{1}\right)  .\medskip$

Let $\Gamma_{0}=S^{n-2}\times S^{1}$ and $p_{0}\in\Gamma_{0}$. Keeping in mind
that $G_{0}=\left\langle a_{0}\right\rangle $, we abuse notation slightly by
letting $a_{0}$ also represent a generator of $\pi_{1}\left(  \Gamma_{0}%
,p_{0}\right)  \cong\mathbb{Z}$. This gives a canonical isomorphism from
$G_{0}$ to $\pi_{1}\left(  \Gamma_{0},p_{0}\right)  $.

Let $\varepsilon$ be a small positive number and $C_{0}^{\prime}=\Gamma
_{0}\times\left[  1-\varepsilon,1\right]  $. To the left-hand boundary
component of $C_{0}^{\prime}$ attach an orientable $1$-handle $h_{0}^{1}$.
Note that $C_{0}^{\prime}\cup h_{0}^{1}$ and its left boundary component each
have fundamental group that is free on two generators--the first corresponding
to $a_{0}$, and the second corresponding to a circle that runs once through
$h_{0}^{1}$. Denote this second generator by $a_{1}$. Keeping in mind the
presentation $G_{1}=\left\langle a_{0},a_{1}\mid a_{1}=[a_{1},a_{0}%
]\right\rangle $, attach to the left-hand boundary component of $C_{0}%
^{\prime}\cup h_{0}^{1}$ a $2$-handle $h_{0}^{2}$ along a regular neighborhood
of a loop corresponding to $a_{1}^{-2}a_{0}^{-1}a_{1}a_{0}$. Let $B_{0}%
=C_{0}^{\prime}\cup h_{0}^{1}\cup h_{0}^{2}$ and let $\Gamma_{1}$ denote the
left-hand boundary component of $B_{0}$. By avoiding the arc $p_{0}%
\times\left[  1-\varepsilon,1\right]  $ when attaching $h_{0}^{1}$ and
$h_{0}^{2}$, we may let $p_{1}=p_{0}\times\left\{  1-\varepsilon\right\}
\in\Gamma_{1}$. Clearly $\pi_{1}\left(  B_{0},p_{1}\right)  \cong G_{1}$. By
inverting the handle decomposition, we may view $B_{0}$ as the result of
attaching an $\left(  n-2\right)  $-handle and then an $\left(  n-1\right)
$-handle to a small product neighborhood $C_{1}$ of $\Gamma_{1}$. Since these
handles have index greater than $2,$ $\Gamma_{1}\hookrightarrow B_{0}$ induces
a $\pi_{1}$-isomorphism. Hence $\pi_{1}\left(  \Gamma_{1},p_{1}\right)  \cong
G_{1} $. This gives us a cobordism $\left(  B_{0},\Gamma_{1},\Gamma
_{0}\right)  $ with the desired boundary components. However, it is not the
cobordism we are seeking.

Next attach a $2$-handle $k_{0}^{2}$ to $B_{0}$ along a circle in $\Gamma_{1}$
representing $a_{1}$. Note that $k_{0}^{2}$ and $h_{0}^{1}$ form a canceling
handle pair in $C_{0}^{\prime}\cup h_{0}^{1}\cup h_{0}^{2}\cup k_{0}^{2}$.
Moreover, since $a_{1}$ has been killed, $h_{0}^{2}$ is now attached along a
trivial loop in the left-hand boundary of $C_{0}^{\prime}\cup h_{0}^{1}\cup
k_{0}^{2}\approx C_{0}^{\prime}$. Provided that $h_{0}^{2}$ was attached with
the appropriate framing (this can still be arranged if necessary), we may
attach a $3$-handle $k_{0}^{3}$ to $C_{0}^{\prime}\cup h_{0}^{1}\cup h_{0}%
^{2}\cup k_{0}^{2}$ that cancels $h_{0}^{2}$. Therefore, $C_{0}^{\prime}\cup
h_{0}^{1}\cup h_{0}^{2}\cup k_{0}^{2}\cup k_{0}^{3}\approx\Gamma_{0}%
\times\left[  0,1\right]  $. The desired cobordism $\left(  A_{0},\Gamma
_{0},\Gamma_{1}\right)  $ will be the complement of $B_{0}$ in this product.
More precisely, $A_{0}=C_{1}\cup k_{0}^{2}\cup k_{0}^{3}$ where $C_{1}$ is a
small product neighborhood of $\Gamma_{1}$ in $B_{0}$. By avoiding $p_{1}$
when attaching $k_{2}$ and $k_{3}$ we may let $p_{0}$ be the left endpoint of
the collar line of $C_{1}$ having right end point corresponding to $p_{1}$. A
schematic diagram of this setup is given in Figure 3.%
\begin{figure}
[ht!]
\begin{center}
\includegraphics[
width=.9\hsize
]%
{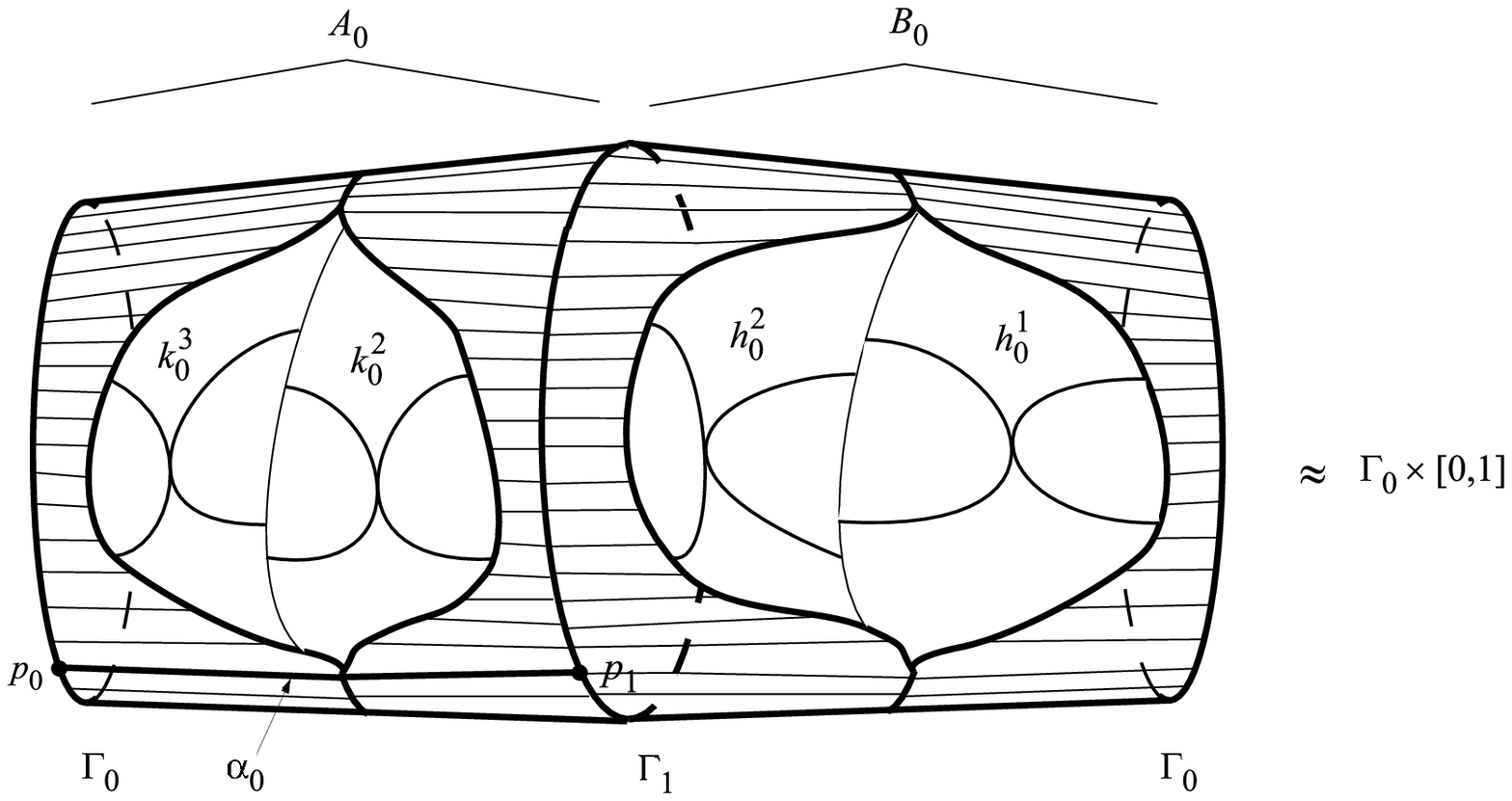}%
\caption{ }%
\end{center}
\end{figure}

By Van Kampen's theorem, it is clear that $\pi_{1}\left(  A_{0},p_{1}\right)
\cong\left\langle a_{0}\right\rangle $, and that the inclusion induced
homomorphism $\pi_{1}(\Gamma_{1},p_{1})\rightarrow\pi_{1}(A_{0},p_{1})$ sends
$a_{0}$ to $a_{0}$ and $a_{1}$ to $1$. By inverting the cobordism, we may view
$A_{0}$ as the result of attaching an $\left(  n-3\right)  $- and an $\left(
n-2\right)  $-handle to the right-hand boundary of $C_{0}=\Gamma_{0}%
\times\left[  0,\varepsilon\right]  $. Hence, inclusion $\Gamma_{0}%
\hookrightarrow A_{0}$ induces the obvious $\pi_{1}$-isomorphism. It follows
that properties a)-c) are satisfied by $\left(  A_{0},\Gamma_{0},\Gamma
_{1}\right)  $.\medskip

\noindent\textbf{Inductive Step}\qua{\sl Construction of }$\left(
A_{j},\Gamma_{j},\Gamma_{j+1}\right)  $.\medskip

Here we assume that $j\geq1$ and that $\left(  A_{j-1},\Gamma_{j-1},\Gamma
_{j}\right)  $ has already been constructed. We will construct $A_{j}$ from
$\Gamma_{j}$ in the same manner that we constructed $A_{0}$ from $\Gamma_{0} $.

Given that $\pi_{1}\left(  \Gamma_{j},p_{j}\right)  \cong G_{j}=\left\langle
a_{0},a_{1},\cdots,a_{j}\mid a_{i}=\left[  a_{i},a_{i-1}\right]  \text{ for
all }1\leq i\leq j\right\rangle $, we expand the fundamental group by
attaching a $1$-handle $h_{j}^{1}$ to the left-hand boundary component of
$C_{j}^{\prime}=\Gamma_{j}\times\left[  1-\varepsilon,1\right]  $. Let
$a_{j+1}$ denote the fundamental group element of $C_{j}^{\prime}\cup
h_{j}^{1}$ corresponding to a loop that runs once through $h_{j}^{1}$. Then
attach to the left-hand boundary component of $C_{j}^{\prime}\cup h_{j}^{1}$ a
$2$-handle $h_{j}^{2}$ along a regular neighborhood of a loop corresponding to
$a_{j+1}^{-2}a_{j}^{-1}a_{j+1}a_{j}$. This yields a cobordism $(B_{j}%
,\Gamma_{j+1},\Gamma_{j})$ with $\pi_{1}\left(  B_{j},p_{j+1}\right)  \cong
G_{j+1}$ and $\Gamma_{j+1}\hookrightarrow B_{j}$ inducing a $\pi_{1}%
$-isomorphism. Now attach a $2$-handle $k_{j}^{2}$ to $B_{j}$ along a circle
in $\Gamma_{j+1}$ representing $a_{j+1}$. Reasoning as in the base case, we
may then attach a $3$-handle $k_{j}^{3}$ to cancel $h_{j}^{2}$ and giving
\[
C_{j}^{\prime}\cup h_{j}^{1}\cup h_{j}^{2}\cup k_{j}^{2}\cup k_{j}^{3}%
\approx\Gamma_{j}\times\left[  0,1\right]  .
\]
Let $C_{j+1}$ be a small product neighborhood of $\Gamma_{j+1}$ in $B_{j}$ and
let
\begin{equation}
A_{j}=C_{j+1}\cup k_{j}^{2}\cup k_{j}^{3}. \tag{\#}%
\end{equation}
Again, the same reasoning used in the base case shows that $\left(
A_{j},\Gamma_{j},\Gamma_{j+1}\right)  $ satisfies conditions a)-c).\medskip

\noindent\textbf{Note}\qua\emph{In completing the proof of Theorem
\ref{perfect}, we will utilize---in addition to properties a-c)---specific
details and notation established in the above construction. }\medskip

It remains to prove the following:

\begin{proposition}
Each clean neighborhood of infinity in $M_{\ast}^{n}$ has finite homotopy type.

\begin{proof}
It suffices to find one cofinal sequence of clean neighborhoods of infinity
with this property. For each $i\geq1$, let $N_{i}^{\prime}=N_{i}\cup
k_{i-1}^{2}$, where $N_{i}=A_{i}\cup A_{i+1}\cup A_{i+2}\cup\cdots$ and
$k_{i-1}^{2}$ is the $2$-handle used in constructing $A_{i-1}$ (See (\#).) We
will show that, for each $i\geq1$, the inclusion
\begin{equation}
\Gamma_{i}\cup k_{i-1}^{2}\hookrightarrow N_{i}^{\prime} \tag{**}%
\end{equation}
is a homotopy equivalence. Hence, $N_{i}^{\prime}$ has finite homotopy type.

Given $i\geq1$, let $A_{i}^{\prime}=A_{i}\cup k_{i-1}^{2}$ and $E_{i}^{\prime
}=A_{i}^{\prime}\cup B_{i}$. Note that $E_{i}^{\prime}$ is not a subset of
$M_{\ast}^{n}$ since $B_{i}$ is not. We now have a cobordism $\left(
E_{i}^{\prime},\Gamma_{i}^{\prime},\Gamma_{i}\right)  $ where (attaching
handles from right to left)
\[
E_{i}^{\prime}=C_{i}^{\prime}\cup h_{i}^{1}\cup h_{i}^{2}\cup k_{i}^{2}\cup
k_{i}^{3}\cup k_{i-1}^{2}\approx\Gamma_{i}\times\left[  0,1\right]  \cup
k_{i-1}^{2}\text{.}%
\]
Here the left-hand boundary $\Gamma_{i}^{\prime}$ may be obtained from
$\Gamma_{i}$ by performing surgery on a regular neighborhood of a circle
representing the element $a_{i}\in\pi_{i}\left(  \Gamma_{i},p_{i}\right)  $.

We may reorder handles so that $k_{i-1}^{2}$ is attached first. (Sliding
$k_{i-1}^{2}$ past $h_{i}^{2}$, $k_{i}^{2}$ and $k_{i}^{3}$ is standard;
attaching $k_{i-1}^{2}$ before $h_{i}^{1}$ requires a quick review of our
construction.) Let $\widehat{k}_{i-1}^{2}\subset int\left(  k_{i-1}%
^{2}\right)  $ be a small regular neighborhood of the core of $k_{i-1}^{2}$,
extended along the product structure of $C_{i}^{\prime}$ to the right-hand
boundary $\Gamma_{i}$. See Figure 4(a).%
\begin{figure}
[ht!]
\begin{center}
\includegraphics[
width=.9\hsize
]%
{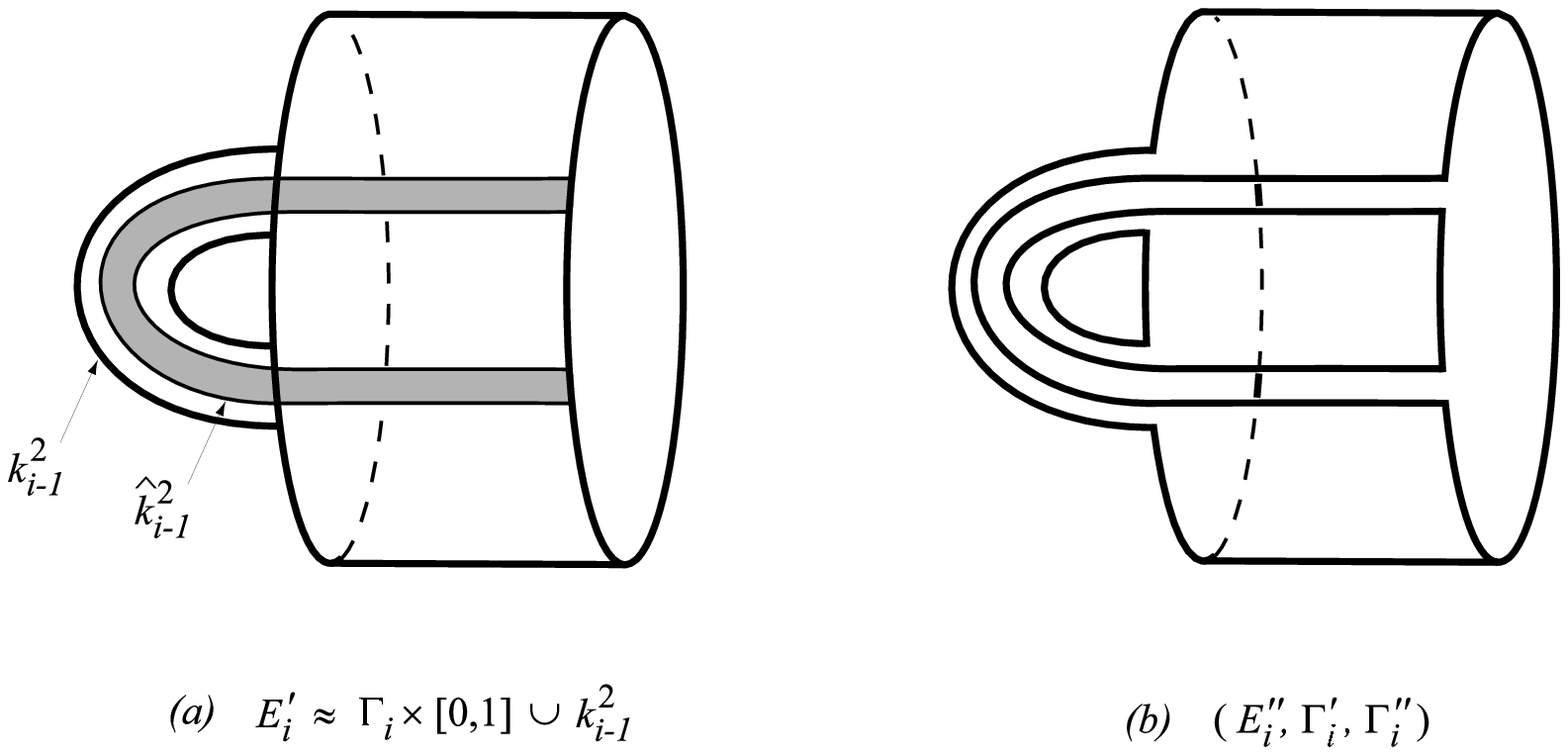}%
\caption{ }%
\end{center}
\end{figure}
Carving from $E_{i}^{\prime}$ the interior of this ``thin'' $2$-handle
$\widehat{k}_{i-1}^{2}$, we obtain a cobordism $\left(  E_{i}^{\prime\prime
},\Gamma_{i}^{\prime},\Gamma_{i}^{\prime\prime}\right)  $ where $\Gamma
_{i}^{\prime\prime}\approx\Gamma_{i}^{\prime}$ since $\Gamma_{i}^{\prime
\prime}$ is obtained from $\Gamma_{i}$ by essentially the same surgery that
produced $\Gamma_{i}^{\prime}$. See Figure 4(b). Furthermore, since $A_{i}\cup
B_{i}\approx\Gamma_{i}\times\left[  0,1\right]  $, it is easy to see that
$E_{i}^{\prime\prime}$ is also a product. The existing handle structure on
$E_{i}^{\prime}$, provides a handle decomposition $E_{i}^{\prime\prime}%
=C_{i}^{\prime\prime}\cup h_{i}^{1}\cup h_{i}^{2}\cup k_{i}^{2}\cup k_{i}^{3}$
where $C_{i}^{\prime\prime}$ is a small product neighborhood of $\Gamma
_{i}^{\prime\prime}$. Recalling that $h_{i}^{2}$ was attached along a circle
in $\Gamma_{i}$ representing $a_{i+1}^{-2}a_{i}^{-1}a_{i+1}a_{i}$ where
$a_{i+1}$ represents a circle that runs once through $h_{i}^{1}$, and noting
that (in $\Gamma_{i}^{\prime\prime}$) $a_{i}$ has been killed by surgery, we
see that $h_{i}^{1}$ and $h_{i}^{2}$ have become a canceling handle pair in
$E_{i}^{\prime\prime}$.

We may split $E_{i}^{\prime\prime}$ as $A_{i}^{\prime\prime}\cup B_{i}%
^{\prime\prime}$ where $B_{i}^{\prime\prime}=C_{i}^{\prime\prime}\cup
h_{i}^{1}\cup h_{i}^{2}$ and $A_{i}^{\prime\prime}$ is obtained from the
left-hand component of $B_{i}^{\prime\prime}$ by attaching $k_{i}^{2}$ and
$k_{i}^{3}$. Alternatively, $A_{i}^{\prime\prime}=A_{i}^{\prime}-int\left(
\widetilde{k}_{i-1}^{2}\right)  $ where $\widetilde{k}_{i-1}^{2}$ is the
interior of a regular neighborhood of the core of $k_{i-1}^{2}$ extended to
the right-hand boundary of $A_{i}^{\prime}$. (The $2$-handle $\widetilde
{k}_{i-1}^{2}$ should be thinner than $k_{i-1}^{2}$, but thicker than
$\widehat{k}_{i-1}^{2}$.) It has already been established that $E_{i}%
^{\prime\prime}$ is a product. Since $h_{i}^{1}$ and $h_{i}^{2}$ form a
canceling pair, $B_{i}^{\prime\prime}$ is also a product. Thus, it follows
from regular neighborhood theory that
\[
A_{i}^{\prime\prime}\approx\Gamma_{i}^{\prime}\times\left[  0,1\right]  .
\]
This last identity will be key to the remainder of the proof.\medskip

\noindent\textbf{Claim}\qua {\sl For each }$i\geq1$\emph{, }$A_{i}^{\prime}%
${\sl\ strong deformation retracts onto }$\Gamma_{i}\cup k_{i-1}^{2}%
${.}\medskip

\textbf{Proof of Claim}\qua It suffices to show that $\Gamma_{i}\cup
k_{i-1}^{2}\hookrightarrow A_{i}^{\prime}$ is a homotopy equivalence. Let
$b_{i-1}^{n-2}$ be a belt disk for $k_{i-1}^{2}$ that intersects the thinner
$2$-handle $\widetilde{k}_{i-1}^{2}$ in a belt disk $\widetilde{b}_{i-1}%
^{n-2}$. By pushing in from the attaching region of $k_{i-1}^{2}$ we may
collapse $\Gamma_{i}\cup k_{i-1}^{2}$ onto $\Gamma_{i}^{\prime}\cup
b_{i-1}^{2}$. See Figure 5.%
\begin{figure}
[ht!]
\begin{center}
\includegraphics[
width=.7\hsize
]%
{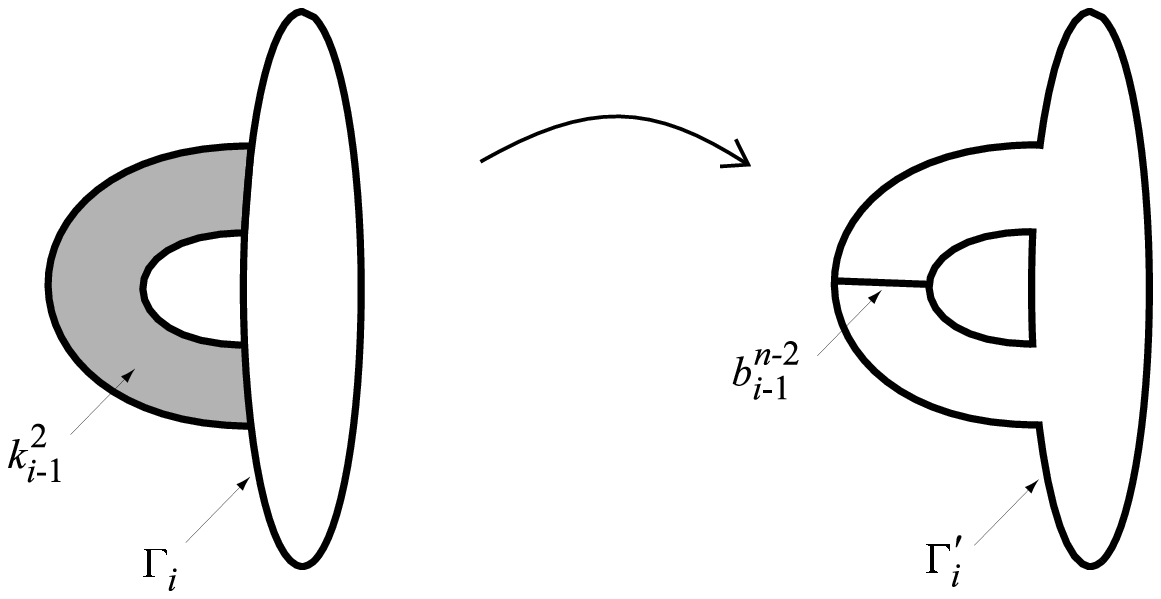}%
\caption{ }%
\end{center}
\end{figure}
Using a similar move, we may collapse $A_{i}^{\prime}$ onto $A_{i}%
^{\prime\prime}\cup\widetilde{b}_{i-1}^{n-2}$. Then, using the product
structure on $A_{i}^{\prime\prime}$ we may collapse $A_{i}^{\prime\prime}%
\cup\widetilde{b}_{i-1}^{n-2}$ onto $\Gamma_{i}^{\prime}\cup b_{i-1}^{n-2}$.
Composing the resulting homotopy equivalences shows that $\Gamma_{i}\cup
k_{i-1}^{2}\hookrightarrow A_{i}^{\prime}$ is a homotopy equivalence and
completes the proof of the claim.\medskip

It is now an easy matter to verify (**). Let $i\geq1$ be fixed. We know that
$A_{i}^{\prime}$ strong deformation retracts onto $\Gamma_{i}\cup k_{i-1}^{2}%
$, and for each $j>i$, we may extend (via the identity) the strong deformation
retraction of $A_{j}^{\prime}$ onto $\Gamma_{j}\cup k_{j-1}^{2}$ to a strong
deformation retraction of $A_{j-1}\cup A_{j}$ onto $A_{j-1}$. By standard
methods, we may assemble these strong deformation retractions to a strong
deformation retraction of $N_{i}^{\prime}$ onto $\Gamma_{i}\cup k_{i-1}^{2}$.
\end{proof}
\end{proposition}

\section{An Open Question}

Our work on pseudo-collars was partly motivated by \cite{CS}, which the
authors advertise as a version of Siebenmann's thesis for Hilbert cube
manifolds. Their result provides necessary and sufficient conditions for a
Hilbert cube manifold $X$ to be ``$\mathcal{Z}$-compactifiable'', ie,
compactifiable to a space $\widehat{X}$ such that $\widehat{X}-X$ is
$\mathcal{Z}$-set in $\widehat{X}$.

\begin{theorem}
[Chapman and Siebenmann]\label{z-compactification}A Hilbert cube manifold $X$
admits a $\mathcal{Z}$-com\-pact\-if\-i\-ca\-tion iff each of the following is satisfied.

\begin{enumerate}
\item[\rm(a)] X is inward tame at infinity.

\item[\rm(b)] $\sigma_{\infty}(X)=0$.

\item[\rm(c)] $\tau_{\infty}\left(  X\right)  \in\underleftarrow{\lim}%
^{1}\left\{  Wh\pi_{1}(X\backslash A)\mid A\subset X\text{ compact}\right\}  $
is zero.
\end{enumerate}
\end{theorem}

Notice that conditions a) and b) are identical to conditions (1) and (3) of
Theorem \ref{pseudo}. The obstruction in c) is an element of the ``first
derived limit'' of the indicated inverse system, where $Wh$ denotes the
Whitehead group functor. See \cite{CS} for details.

It is not well-understood when conditions a)-c) imply $\mathcal{Z}%
$-compactifiability for spaces that are not Hilbert cube manifolds. In
\cite{Gu2}, a polyhedron was constructed which satisfies the hypotheses of
Theorem \ref{z-compactification}, but which fails to be $\mathcal{Z}%
$-compactifiable. However, it is unknown whether a finite dimensional manifold
that satisfies these conditions can always be $\mathcal{Z}$-compactified. In
trying to answer this question, it seems worth noting that Chapman and
Siebenmann employed a two step procedure in proving their result. First they
showed that a Hilbert cube manifold satisfying conditions a) and b) is
pseudo-collarable. Next they used the pseudo-collar structure, along with
condition c) and some powerful Hilbert cube manifold techniques to obtain a
$\mathcal{Z}$-compactification.

In contrast with the infinite dimensional situation, the manifolds $M_{\ast
}^{n}$ constructed in this paper satisfy conditions a) and b) yet fail to be
pseudo-collarable. Furthermore, an inductive application of the exact sequence
on page 157 of \cite{Wld} shows that each group $G_{i}$ appearing in the
canonical inverse sequence representative of $\pi_{1}\left(  \varepsilon
(M_{\ast}^{n})\right)  $ has trivial Whitehead group. It follows that
$\tau_{\infty}\left(  M_{\ast}^{n}\right)  =0$. Thus, the $M_{\ast}^{n}$'s
would appear to be ideal candidates for counterexamples to an extension of
Theorem \ref{z-compactification} to the case of finite dimensional manifolds.
More generally, we ask:\smallskip

\noindent\textbf{Question}\ {\sl Can a }$\mathcal{Z}${\sl-compactifiable
open }$n${\sl-manifold fail to be pseudo-collarable?}

\end{document}